\documentstyle[11pt,psfig]{article}
\hsize \textwidth      
\advance \hsize by -\marginparwidth
\evensidemargin \oddsidemargin  

\advance\hoffset by 5mm 

\newtheorem{thm}{Theorem}
\newtheorem{prop}{Proposition}
\newtheorem{lemma}{Lemma}
\newtheorem{corollary}{Corollary}
\newcommand{\Rm}{I\!\!R}
\newcommand{\pdr}[2]{\frac{\partial{#1}}{\partial{#2}}}
\begin{document}
\title{Enhancement of the Traveling Front Speeds in 
 Reaction-Diffusion Equations with Advection}

\author{Alexander Kiselev  and Leonid Ryzhik\\
Department of Mathematics \\ University of Chicago\\
Chicago IL 60637}

\maketitle
\begin{abstract}
  We establish rigorous lower bounds on the speed of traveling fronts
  and on the bulk burning rate in reaction-diffusion equation with passive
   advection. The non-linearity is assumed to be of
  either KPP or ignition type. We consider two main classes of flows.
  Percolating flows, which are characterized by the presence of long
  tubes of streamlines mixing hot and cold material, lead to strong speed-up
  of burning which is linear in the amplitude of the flow, $U$.
  On the other hand the cellular flows, which have closed streamlines,
  are shown to produce weaker increase in reaction.  For such flows
  we get a lower bound which grows as $U^{1/5}$ for a large
  amplitude of the flow.
\end{abstract}

\section{Introduction}
Propagation of thin fronts in moving fluids arises in many 
situations in physics and engineering.
Consider a mixture of reactants interacting in a region that may have
 a rather complicated spatial structure but is thin across.
The reaction front moves towards the unburned reactants leaving behind the
burned ones. When the reactants are mixed by an ambient fluid
then the burning rate may be enhanced. The physical reason
for this observed speed-up is believed to be that fluid
advection tends to increase the area available for reaction.
Many important engineering applications of
combustion operate in the presence of turbulent advection,
and therefore the influence of advection on burning has been
studied extensively by physicists, engineers and mathematicians.
In the physical literature one can find a number of models and
approaches that yield different predictions --
relations between the turbulent intensity and the burning rate
\cite{CW,KA,ker-1,Y}. These results are usually obtained using
heuristic models and physical reasoning. For a recent review of some
of the physics literature we refer to \cite{curved-fronts,Ro}.

The lack of agreement between different physical models makes
rigorous results, even for simplified mathematical models, particularly
valuable and useful.
A well-established mathematical model
that describes a chemical reaction in a fluid is a system of two
equations for concentration $C$ and temperature $T$ of the form
\begin{eqnarray}\label{gen-sys}
 && T_t+u\cdot\nabla T=\kappa\Delta T+\frac{v_0^2}{\kappa}g(T)C\\
&&  C_t+u\cdot\nabla C=\frac{\kappa}{\hbox{Le}}\Delta C-
\frac{v_0^2}{\kappa}g(T)C.\nonumber
\end{eqnarray}
For exposition purposes, all consideration in this paper will be
carried out in two spacial dimensions, but our methods extend to an
arbitrary dimension in a straightforward way. Equations
(\ref{gen-sys}) are coupled to the reactive Euler equations for the
advection velocity $u(x,y,t)$.  Two assumptions are usually made to
simplify the problem: first, constant density approximation \cite{CW}
that allows to decouple the Euler equations from the system
(\ref{gen-sys}). Then one may consider $u(x,y,t)$ as a prescribed
quantity that does not depend on $T$ and $C$.  Furthermore, it is
often assumed that $\hbox{Le}=1$, or, equivalently, thermal and
material diffusivities are equal. These two assumptions allow to
reduce the above system to a single scalar equation for the
temperature $T$:
\begin{equation}
  \label{eq:0.1}
  \frac{\partial T}{\partial t} + u(x,y,t)\cdot\nabla T=\kappa\Delta T+
\frac{v_0^2}{\kappa}f(T)
\end{equation}
with $f(T)=g(T)(1-T)$, provided that $C(x,y,0)=1-T(x,y,0)$.  We will
consider the problem (\ref{eq:0.1}) in a strip
$\Omega=\Rm_x\times[0,H]_y$ with boundary conditions in $x$:
\begin{equation}
  \label{eq:0.2}
  T(x,y,t)\to 1 \hbox{ as $x\to -\infty$},~~~~
T(x,y,t)\to 0 \hbox{ as $x\to +\infty$} 
\end{equation}
and either Neumann 
\begin{equation}
  \label{eq:0.3}
  \pdr{T}{y}(x,0,t)=\pdr{T}{y}(x,H,t)=0
\end{equation}
or periodic
\begin{equation}
  \label{eq:0.4}
  T(x,y,t)=T(x,y+H,t)
\end{equation}
boundary conditions in $y$. Furthermore, we assume that the initial
data $T_0(x,y)$ for (\ref{eq:0.1}) satisfies the bounds
\begin{eqnarray}
  \label{eq:in-decay}
 && T_0(x,y)=1-O(e^{\lambda x})~~\hbox{for}~~x<0,~~
  T_0(x,y)=O(e^{-\lambda x})~~\hbox{for}~~x>0, \\
\label{eq:1.3a} &&|\nabla T_0|= O(e^{-\lambda|x|})~~
\hbox{for}~~\hbox{some}~~\lambda>0.
\end{eqnarray}

We adapt fairly general assumptions on
$f,$ requiring only that $f(T)$ is not equal identically to zero, and
\begin{eqnarray}
  \label{eq:nonlin}
  f(0)=f(1)=0,~~f(T)\ge 0~~\hbox{ for $T\in(0,1)$},~~f\in C^1[0,1].
\end{eqnarray}
Two types of reaction rates $f(T)$ are
distinguished in this class. The KPP-type reactions satisfy 
\begin{equation}\label{eq:kpp}
f(0)=f(1)=0,~~\hbox{$f(T)>0$ for
$T\in(0,1).$}
\end{equation}
An additional requirement $f'(0)={\rm max}_{T \in [0,1]}f(T)/T$
is often made. We do not make such requirement in this paper,
and call the class described by (\ref{eq:kpp}) general KPP.
Our interpretation of KPP includes an important Arrhenius-type
non-linearity, 
\[ 
f(T) = C(1-T)e^{-A/T}, 
\]
that is believed to be an appropriate model for many 
chemical reactions in the context of reaction-diffusion models. 
We also consider the ignition non-linearities with  
\begin{equation}\label{eq:ignit}
f(T)=0~~\hbox{ for
$T\in[0,\theta_0]$ and $T=1$,~ $f(T)>0$ for $T\in(\theta_0,1)$.}
\end{equation}
By our assumptions on the nonlinearity (\ref{eq:nonlin}), we
can find $\theta_4 >\theta_1,$ and $f_0,\zeta>0$ such that
\begin{equation}
  \label{eq:cond1}
  f(\theta)>f_0~~\hbox{for $\theta\in(\theta_1-\zeta,\theta_4+\zeta)$}.
\end{equation}
 The values of the constants $f_0$,
$\zeta$ and $\theta_{1,4}$ are  
the only information on the nonlinearity $f(T)$ that shows up in our
bounds on the burning rate.

We assume that advection $u(x,y)\in
C^1(\Omega)$ is time independent, has mean zero in the
$x$-direction:
\begin{equation}
  \label{eq:mean-zero}
  \int\limits_{0}^Hu_1(x,y)dy=0
\end{equation}
and is incompressible:
\begin{equation}
  \label{eq:incomp}
  \nabla\cdot u=0.
\end{equation}

The mathematical literature on the scalar reaction-diffusion equation
(\ref{eq:0.1}) is enormous; far from giving an exhaustive overview, we
mention several papers directly related to our work.  First rigorous
results about traveling waves for equation (\ref{eq:0.1}) go back to
classical works of Kolmogorov, Petrovskii and Piskunov \cite{KPP} and
Fisher \cite{Fisher}, which considered the case $u=0$ in one dimension
for the KPP nonlinearity.  Recently equation (\ref{eq:0.1}) with $u\ne
0$, and in particular the effect of advection, became a subject of
intense research.  Berestycki and Nirenberg \cite{Ber-Nir-1,Ber-Nir-2}, 
and Berestycki, Larrouturou and Lions
\cite{Ber-Lar-Lions} initiated the studies of the existence of traveling
waves for equation (\ref{eq:0.1}) of the form
\begin{equation}
  \label{eq:shear-front}
  T(x,y,t)=T(x-ct,y),
\end{equation}
for shear flows of the form $u=(u(y),0)$.  Their stability was studied in
\cite{Ber-Lar-Roq,Mal-Roq,Roq-1}, while in further works
\cite{Ber-Hamel,jxin-1,jxin-2} stability
and existence of traveling waves were established
for the wider class of periodic flows. In this case,
the traveling fronts have the form
 \[  T(x,y,t)=U(x-ct,x,y) \]
and are periodic in the last two variables.  These and other results
were recently reviewed in \cite{jxin-3}, and we refer the reader to
this paper for a detailed exposition of the subject.  Until very
recently, there were no rigorous results on the physically interesting
question of the speed of traveling waves.  First such results have
been established in \cite{CKOR} for percolating flows, and in
\cite{ABP} and \cite{HPS} for the shear flows.  Numerical studies
of the propagation of fronts were performed for a shear flow in
\cite{siv-kogan-1} with ${\hbox{Le}}\ne 1$, and for cellular flows in
\cite{siv-kogan-2}.

Another major direction of research has been homogenization approach.
The homogenization regime $\kappa \to 0,$ when the front width goes to
zero, was extensively studied for KPP-type nonlinearity and for
advection velocity that is periodic and varies either on the integral
or diffusive scale by Freidlin \cite{freidlin-1, freidlin-2,
  freidlin-3}.  Recently Majda and Souganidis derived an effective
Hamilton-Jacobi equation in the limit $\kappa \to 0$ for the case of
advection velocity varying on a small $\kappa-$dependent scale that is
larger or comparable to that of the front width \cite{MS}.  This
effective equation is still difficult to analyze, and analytical
predictions have been derived only for the shear flows. Numeric
experiments exploring the results of \cite{MS} have been carried out
in \cite{Emb-Maj-Sou-1,Emb-Maj-Sou-2,McL-Zhu}.

Very recently, Hamel \cite{Hamel} and Heinze, Stevens and Papanicolaou
\cite{HPS} proposed an elegant variational approach to the estimates 
of the speed of traveling waves in the presence of periodic advection.
However, to the best
of our knowledge nontrivial lower bounds using this method were
obtained so far only for shear flows in the homogenization regime or
for small advection, where they provide precise bounds for the 
small speed-up of the front \cite{HPS}.

The key question we wish to address in this paper 
is: what characteristics of the
ambient fluid flow are responsible for burning rate enhancement?
The question needs first to be made precise, because the reaction region
may be complicated and, in general, may move with an ill-defined
velocity, when traveling fronts do not exist. 
To measure the speed of burning in such situations,  the bulk burning rate
\begin{equation}
  \label{eq:0.5}
  V(t)=\int\limits_{\Omega}T_t(x,y,t)\frac{dxdy}{H}
\end{equation}
and its  time average
\[  \langle V\rangle_t=\frac{1}{t}\int\limits_0^tV(s)ds \]
have been recently introduced in \cite{CKOR}.
Note that for traveling fronts of the form (\ref{eq:shear-front}) we
have $V(t)=c$, but the notion of bulk burning rate 
makes sense in much more general
situations when traveling fronts of the above form may not exist, and
bulk burning rate serves as a natural generalization of the front
speed. We have obtained in \cite{CKOR} lower bounds for
$
\displaystyle\langle V\rangle_t
$
when $f(T)$ is a concave function of the KPP type. The bounds are
linear in the magnitude of the advecting velocity $u(x,y)$ provided
that there exist tubes of streamlines that connect $x=-\infty$ and
$x=+\infty,$ satisfying some mild additional technical assumptions.
We say then that the flow is percolating. In particular
these bounds hold for shear flows of the form $(u(y),0)$. 

In this paper we consider much more general reaction rates $f(T)$ that
are either of the ignition or general KPP type, and establish similar
lower bounds for $V(t)$ for percolating flows that are periodic in
space. The bound is linear in the magnitude of $u$ and deteriorates as
the scale of oscillations of $u$ becomes comparable to the laminar
front width $l=\kappa/v_0.$ 
It is easy to show \cite{CKOR} that for any $u \in C^1,$ the burning
rate $\langle V \rangle_t$ satisfies linear in $\|u\|_\infty$ 
upper bound (for initial data as in (\ref{eq:in-decay}) and (\ref{eq:1.3a})).
Therefore, shear (and, more generally, percolating) flows are
as effective as possible in speeding up combustion in terms of 
the power of $\|u\|_\infty$ in the large intensity regime.  
\begin{figure}
  \centerline{
  \psfig{figure=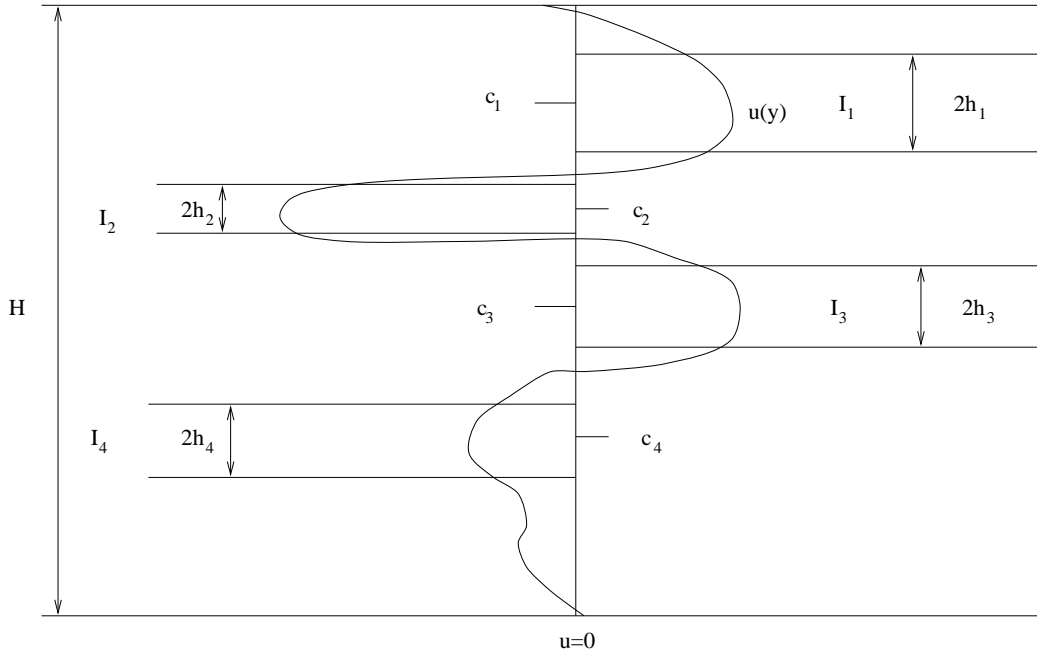,width=16cm} 
     }
  \caption{The structure of the shear flow}
  \label{fig0}
\end{figure}
In particular, we establish the following lower bound for
the bulk burning rate in a shear flow.
\begin{thm}\label{thm0.1}
  Let $T_0(x,y)$ be an arbitrary initial data satisfying
  (\ref{eq:in-decay}) and (\ref{eq:1.3a}), and let $T(x,y,t)$ satisfy
  (\ref{eq:0.1}) with the either the Neumann (\ref{eq:0.3}) or
  periodic boundary conditions (\ref{eq:0.4}). Let also
  $u(x,y)=(u(y),0)$ in (\ref{eq:0.1}). Then both for KPP and ignition
  non-linearities we have
  \begin{equation}
    \label{eq:cor1.1}
    \lim_{t\to\infty}\langle V\rangle_t\ge 
C\left(v_0+\sum_j\left(1+\frac{l}{h_j}\right)^{-1}
\int\limits_{c_j-h_j}^{c_j+h_j}|u(y)|\frac{dy}{H}\right), 
  \end{equation}
where the constant  $C$ depends only on the reaction function $f$. Here
the intervals $I_j=[c_j-h_j,c_j+h_j]\in[0,H]$ are any intervals such that
\begin{equation}\label{eq:def-Ij}
\frac{\|u\|_{\infty,j}}{2}\le |u(y)|\le {\|u\|_{\infty,j}},~~~
{\|u\|_{\infty,j}}=\sup_{I_j}|u(y)|.
\end{equation}
We do not require $\cup_j I_j = [0,H].$ 
\end{thm}
The choice of intervals $I_j$ is up to us, and should be 
made to maximize the lower bound.
See Figure~\ref{fig0} for an illustration. 

As a corollary, the bound (\ref{eq:cor1.1}) holds for the speed $c$ of
a traveling front of the from (\ref{eq:shear-front}).  Our bound
behaves correctly in the homogenization regime when $u(y)$ has the
form $u(y)=\frac{A}{\varepsilon}v(y/\varepsilon)$ and provides a bound
that is linear in the magnitude $A$ of advection, in agreement with
\cite{CKOR,HPS}, where homogenization limit was studied. We also prove
the analog of Theorem \ref{thm0.1} for general percolating flows (see
Theorem \ref{cor3} in Section~\ref{sec:percol}). 

Another main result of this work concerns cellular flows with closed
streamlines.  Roughly speaking, in terms of their burning enhancement
properties, such flows can be thought of as ``the worst'' class of
flows, opposing ``the best'' percolating flows.  One can expect the
burning enhancement to be significantly weaker for cellular flows
because of the numerous diffusive interfaces which prevent hot and
cold regions from mixing fast.  Cellular flows pose mathematically
more challenging problem because of these diffusive interfaces; we
will see that the estimates for percolating flows will form only a
fraction of the argument we will need in the cellular case. We
consider a particular example of a cellular flow
\begin{figure}
  \centerline{
  \psfig{figure=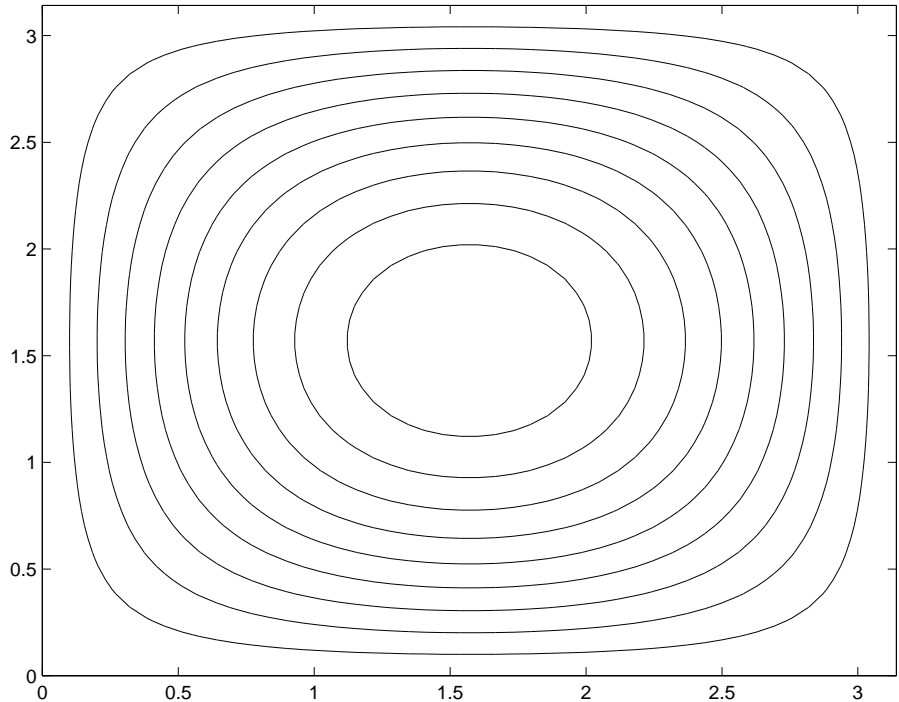,width=12cm} 
     }
  \caption{Streamlines of a cellular flow ($H=\pi$).}
  \label{fig:cell}
\end{figure}
\begin{equation}
  \label{stream0}
  u(x,y)=UH\left(\pdr{\psi}{y},-\pdr{\psi}{x}\right),~~
\psi(x,y)=\sin(\frac{x}{H})\sin(\frac{y}{H}).
\end{equation}
See Figure~\ref{fig:cell} for an illustration of streamlines in a
single cell.
Then under the assumptions of a large Peclet number
and thin laminar front width:
\begin{equation}\label{eq:assump0}
\hbox{Pe}=\frac{UH}{\kappa}\ge 1,~~~
  \frac{l}{H}\le 1, ~~
l=\frac{\kappa}{v_0}
\end{equation}
we get the following result.
\begin{thm}
\label{celltravelthm} 
Let $T_0(x,y)$ be an arbitrary initial data satisfying
(\ref{eq:in-decay}) and (\ref{eq:1.3a}), and let $T(x,y,t)$ satisfy
(\ref{eq:0.1}) with the either the Neumann (\ref{eq:0.3}) or periodic
boundary conditions (\ref{eq:0.4}). Let also $u(x,y)$ in
(\ref{eq:0.1}) be given by (\ref{stream0}), and assume (\ref{eq:assump0})
is satisfied. Then both for KPP and
ignition non-linearities we have
\begin{equation}
\lim_{t\to\infty}\langle V\rangle_t\ge\left\{
{\begin{matrix}{
\left(C_1\sqrt{\frac{\tau_c}{\tau_u}}+C_2\right)v_0, &
{\rm if}~~  
{\tau_c}\le\tau_u\cr
\left(C_1\left(\frac{\tau_c}{\tau_u}\right)^{1/5}+C_2\right)v_0, &
{\rm if}~~  
{\tau_c}\ge\tau_u\cr}\end{matrix}}\right.
\end{equation}
Here $\tau_c=\kappa/v_0^2$ is the chemical reaction time and
$\tau_u=H/U$ is the turnover time.  The constants in the inequalities
depend only on the reaction $f$, more particularly on the constants $f_0,$
$\zeta,$ and $\theta_4-\theta_1$ that appear in (\ref{eq:cond1}).
\end{thm}
To the best of our knowledge this is the first rigorous bound on the
traveling front speed in a cellular flow. Note that the change of behavior of
our bound depending on the ratio $\tau_c/\tau_u$ is physically natural
since for $\tau_u\ll\tau_c$ the front folds onto itself inside the
period cell, which diminishes the affect of advection.  The lower
bound of Theorem~\ref{celltravelthm} displays square root dependence
on the flow intensity $U$ until $U$ reaches a critical value
determined by a condition $\tau_c = \tau_u.$ After that, the lower
bound behaves like $U^{1/5}.$ Recently, Audoly, Berestycki and Pomeau
\cite{ABP} gave an heuristic argument which proposes that the speed of
the traveling front for cellular flows should scale as $U^{1/4}$ in
the large $U$ limit, which may indicate that our lower bound is not
far off from the sharp bound.

One of the fundamental mathematical difficulties we deal with in this paper
may be roughly described as follows. We will be able to bound the burning rate
from below by integrals over the domain of reaction term $f(T)$ and of the 
square of the gradient $|\nabla T|^2.$ It turns out that in order to
obtain a lower bound on $V$ in terms of $u,$ it will be necessary
to bound the integral of the higher derivative Laplacian term in terms of 
integrals of $f(T)$ and $|\nabla T|^2.$ One can expect to do this 
using parabolic regularity, but the constants in such a priori estimates
typically depend on $u$, and this dependence turns out to be too crude to get 
interesting results. We tackle this difficulty by taking advantage of
the fact that what we need to estimate is the integral of Laplacian,
not of the absolute value of Laplacian, and employ an appropriate averaging
procedure to reduce derivatives. We hope that this idea will be useful 
in other related contexts in PDE estimates. 

The paper is organized as follows. We prove Theorem \ref{thm0.1} in
Section \ref{sec:shear}, as well as some other results for shear flows.
The analogous results for the percolating flows are proved in Section
\ref{sec:percol}. We present our main results for the cellular flows,
in particular implying Theorem \ref{celltravelthm}, 
in Section \ref{sec:cell-main}.
Sections \ref{sec:advect}, \ref{sec:diff} and \ref{sec:react} contain
some intermediate estimates in the proof of Theorem \ref{celltravelthm}. We
put these estimates together to finish the proof in Section
\ref{sec:putting}. 

\noindent {\bf Acknowledgment} \\
We wish to thank Peter Constantin for initiating our collaboration on
the subject, and for numerous useful discussions and advice.  We are
grateful to Henri Berestycki, Fausto Cattaneo, Albert Fannjiang,
Francois Hamel and George Papanicolaou for interesting discussions.
AK was supported by NSF grant DMS-9801530, and LR was supported by NSF
grant DMS-9971742 and ASCI Flash Center at the University of Chicago.

\section{Shear flows and general nonlinearities}\label{sec:shear}

We first consider (\ref{eq:0.1}) in a shear, or unidirectional, flow
$(u(y),0)$, which is a particular example of a percolating flow. The
proofs are somewhat less technical in this case and allow us to
introduce some of the ideas used in the general case.  Equation
(\ref{eq:0.1}) in a shear flow becomes
\begin{eqnarray}
  \label{eq:1.1}
&&  \pdr{T}{t}+u(y)\pdr{T}{x}=\kappa\Delta T+\frac{v_0^2}{\kappa}f(T)\\
&&T(x,y,0)=T_0(x,y).\nonumber
\end{eqnarray}
The  advection is assumed to be mean-zero:
\begin{equation}
  \label{eq:shear-mean-zero}
  \int\limits_{0}^Hu(y)dy=0. 
\end{equation}
We impose an additional assumption 
\begin{equation}
  \label{eq:1.2}
  \pdr{T}{t}\ge 0.
\end{equation}
This condition is satisfied for all times 
provided that initially we have
\begin{equation}
  \label{eq:1.2.1}
  \kappa\Delta T_0-u(y)\pdr{T_0}{x}+\frac{v_0^2}{\kappa}f(T_0)\ge 0,
\end{equation}
as follows from the maximum principle (see, e.g. \cite{PW}). 
Therefore (\ref{eq:1.2}) is not
a constraint on the dynamics but rather on the initial data. We note
that (\ref{eq:1.2}) is true for traveling fronts of the form
$T(x-ct,y)$ both for general KPP and ignition nonlinearity
\cite{Ber-Lar-Lions,Ber-Nir-2}. We assume the usual boundary conditions
(\ref{eq:0.2}) at the left and right ends of the strip $\Omega$ and
either Neumann (\ref{eq:0.3}) or periodic (\ref{eq:0.4}) boundary
conditions at $y=0,H$. We also require that the initial data
$T_0(x,y)$ satisfies
(\ref{eq:in-decay}) and (\ref{eq:1.3a}).
These conditions are preserved by evolution (see e.g. \cite{CKOR})
if the advecting velocity $u(y)\in C^1[0,H]$, that is, we have for each $t>0$:
\begin{eqnarray}
  \label{eq:later-decay}
&&  1-T(x,y,t)\le C(t)e^{\lambda x}~~\hbox{for}~~x<0,~~
  T(x,y,t)\le C(t)e^{-\lambda x}~~\hbox{for}~~x>0, \\
&& |\nabla T(x,y,t)|\le C(t)e^{-\lambda|x|}~~\nonumber
\end{eqnarray}  
provided that (\ref{eq:in-decay}) and (\ref{eq:1.3a}) hold initially.

Let $I_j=(c_j-h_j,c_j+h_j)\subset [0,H]$ be a collection of intervals
satisfying (\ref{eq:def-Ij}).
In particular $u(y)$ does not change sign on the intervals $I_j$.
We do {\it not} require that $\displaystyle \cup_jI_j=[0,H]$.
Then the bulk burning rate $V(t)$ defined by (\ref{eq:0.5}) obeys a lower bound
described
by the following Theorem, which is the first main result of this section.
\begin{thm}\label{thm1} 
  Let $T(x,y,t)$ be a solution of (\ref{eq:1.1}) with the boundary
  conditions (\ref{eq:0.2}) and either (\ref{eq:0.3}) or
  (\ref{eq:0.4}). Let the initial data $T_0(x,y)$ satisfy
  (\ref{eq:1.2.1}), (\ref{eq:in-decay}) and (\ref{eq:1.3a}).
  Furthermore, assume that $u(y)\in C^1[0,H]$ has mean zero
  (\ref{eq:shear-mean-zero}) and the nonlinearity $f(T)$ satisfies
  (\ref{eq:nonlin}).  Then there exists a constant $C>0$ that depends
  on $f(T)$ but not on $T_0(x,y)$ or $u(y)$, such that for any
  collection of intervals $I_j$ that satisfies (\ref{eq:def-Ij}) we
  have
  \begin{eqnarray}
    \label{eq:sh-bound}
    V(t)\ge C\left( v_0 +\sum_j\left(1+\frac{l}{h_j}\right)^{-1}
\int\limits_{c_j-h_j}^{c_j+h_j}|u(y)|\frac{dy}{H} \right)
\end{eqnarray}
with $\displaystyle l=\frac{\kappa}{v_0}$.
\end{thm}
{\it Remark.} 1. The lower bound (\ref{eq:sh-bound}) does not deteriorate
when oscillations in $u(y)$ become faster in space as long as its amplitude
grows according to
\[
\frac{\|u\|_{\infty,j}}{v_0}=O\left(\frac{l}{h_j}\right).
\]
This agrees well with the homogenization limit $\displaystyle
u_\varepsilon (y)=\frac 1\varepsilon u({y}/{\varepsilon})$ considered
in \cite{CKOR} and \cite{HPS}, that produces speed-up of the front of
order $O(v_0)$. This is also an improvement of the analogous lower bound
for $V(t)$ for the convex KPP case obtained in \cite{CKOR}, where
$(l/h_j)^2$ appeared in the factor.

2. The regularity assumption on $u(y)$ is used only to guarantee
preservation of the boundary conditions (\ref{eq:1.3a}) that allows us
to integrate by parts in the proof. None of our bounds depend on the
size of derivatives of $u(y)$.

It has been shown in \cite{Ber-Lar-Lions,Ber-Nir-2} both in the case
of ignition non-linearity (\ref{eq:ignit}), and for the general KPP
nonlinearity (\ref{eq:kpp}) that there exist traveling front solutions
of (\ref{eq:1.1}) of the form $T(x,y,t)=U(x-ct,y)$.  The speed $c=c_*$
is uniquely determined by the nonlinearity $f(T)$ and advection $u(y)$
in the ignition case, while traveling front solutions exist for $c\ge
c_u$ for some minimal speed $c_u$ in the KPP case. The function $U(s,y)$ is
monotonically decreasing in the variable $s=x-ct$ in both cases, so
that (\ref{eq:1.2}) holds.
Theorem \ref{thm1} implies the following estimate on the speeds of the
traveling fronts.
\begin{corollary}\label{cor1}
  Let $T(x-ct,y)$ be a traveling front solution of (\ref{eq:1.1})
  with $f(T)$ being either of the ignition nonlinearity type
  (\ref{eq:ignit}), or of the KPP type (\ref{eq:kpp}). Then there exists a
  constant $C>0$ that depends on the function $f$ but not on $u(y)$
  such that
\begin{equation}
  \label{eq:tr-wave-bd}
  c\ge C\left( v_0 +\sum_j\left(1+\frac{l}{h_j}\right)^{-1}
\int\limits_{c_j-h_j}^{c_j+h_j}|u(y)|\frac{dy}{H} \right)
\end{equation}
\end{corollary}
Corollary \ref{cor1} follows immediately from Theorem \ref{thm1} since
we have $V(t)=c$ for $T(x,y,t)=U(x-ct,y)$ due to the boundary
conditions $U(s,y)\to 1$ as $s\to -\infty$, $U(s,y)\to 0$ as $s\to
+\infty$. Corollary 1 and the stability results for traveling fronts
\cite{Roq-1,jxin-2} 
 imply Theorem \ref{thm0.1} for general initial data.
We prove now Theorem \ref{thm0.1} assuming the result of Corollary \ref{cor1}.
\newline
{\bf Proof of Theorem \ref{thm0.1}.} Consider first ignition non-linearity. 
In this case we will show that
\[
\lim_{t\to\infty}\langle V\rangle_t=c_*,
\]
where $c_*$ is the unique speed of the traveling front.
Then (\ref{eq:cor1.1}) will follow from (\ref{eq:tr-wave-bd}).
It was shown in \cite{jxin-2} that for the initial data
satisfying (\ref{eq:in-decay}) (actually just tending to $1$ and $0$
at the two ends) there exist functions $\xi_{1,2}(t)$, 
such that
\[
|\xi_i(t)|=o(t),~~as~~ t\to\infty,
\]
and functions $q_i(t,x)$ that satisfy the linearized problem
\[
\pdr{q_i}{t}+u(y)\pdr{q_i}{x}=\Delta q_i
\]
such that
\begin{equation}
\label{asymvel}
U(x-c_*t+\xi_1(t),y)-q_1(x,y,t)\le T(t,x,y)\le U(x-c_*t-\xi_2(t),y)+q_2(x,y,t).
\end{equation}
Here $U(x-c_*t,y)$ is the traveling wave solution of (\ref{eq:1.1}).
The initial data $q_i(x,y,0)$ may be chosen in $L^1\cap L^\infty(\Omega)$
Then we have for any $c>c_*$:
\begin{eqnarray*}
&&  \langle V\rangle_\tau=\frac 1\tau\int\limits_0^\tau dt\int\limits_{\Omega}
\frac{dxdy}{H}
T_tdxdy=\frac 1\tau\int\limits_{\Omega}\frac{dxdy}{H}[T(x,y,\tau)-T_0(x,y)]\\
&&=
\frac 1\tau\int\limits_{-\infty}^0dx\int\limits_0^H\frac{dy}{H}[
(1-T_0)-(1-T)]+
\frac 1\tau\int\limits_{0}^{c\tau}dx\int\limits_0^H\frac{dy}{H}[T(x,y,\tau)-T_0(x,y)]\\
&&+
\frac 1\tau\int\limits_{c\tau}^{\infty}dx\int\limits_0^H\frac{dy}{H}[T(x,y,\tau)-T_0(x,y)]
\\
&&\le
\frac{C}{\tau}+c+
\frac 1\tau\int\limits_{(c-c_*)\tau}^{\infty}dx\int\limits_0^H\frac{dy}{H}
U(x-\xi_2(t),y)
+\frac 1\tau\int\limits_{c\tau}^{\infty}dx\int\limits_0^H\frac{dy}{H}q_2(x+c_*t,y,t)\le 
\frac{C'}{\tau}+c
\end{eqnarray*}
and hence $\displaystyle\limsup_{\tau\to\infty}\langle
V\rangle_\tau\le c_*$.  Similarly one may show that
$\displaystyle\liminf_{\tau\to\infty}\langle V\rangle_\tau\ge c'$ for
any $c'<c_*$, which shows that (\ref{eq:cor1.1}) holds.

In the KPP case, the estimates of the sort (\ref{asymvel}) are not yet
available. However the bound (\ref{eq:cor1.1}) can be shown by
reduction to the ignition non-linearity case. Indeed, given KPP type
reaction $f,$ consider ignition type reaction $f_\theta=\chi_\theta f
\leq f,$ say by cutting $f$ off in a small neighborhood near zero for
$T\le\theta$. The constant $C$ in (\ref{eq:sh-bound}) does not depend
on $\theta$ for $\theta$ small enough as will be seen from the proof of
Theorem \ref{thm1}. Let $T$ and $T_\theta$ satisfy equations with
reactions $f$ and $f_\theta$ respectively with the same initial data
$T_0(x,y)$. Then $Z = T-T_\theta$ satisfies
\[ 
Z_t +u\cdot\nabla Z -\kappa \Delta Z = \frac{v_0^2}{\kappa}(f(T)-f_\theta
(T_\theta))\ge \frac{v_0^2}{\kappa}(f(T)-f(T_\theta)). 
\]
It follows from the maximum principle that if $Z(x,0) \geq 0$ then
$Z(x,t) \geq 0.$ Hence for the same initial data, the burning rate for
the KPP reaction $f$ is not smaller than for the ignition
non-linearity $f_\theta$:
\[
\langle V[T]\rangle_t=\frac 1t\int \frac{dxdy}{H}[T(x,y,t)-T_0(x,y)]
\ge\frac 1t\int \frac{dxdy}{H}[T_\theta(x,y,t)-T_0(x,y)]
=\langle V[T_\theta]\rangle_t.
\] 
 This implies the validity of the lower bound
(\ref{eq:cor1.1}) in the KPP case.  $\Box$

We now turn to the proof of Theorem \ref{thm1}.  
The proof follows the general ideas of
\cite{CKOR} with significant modifications required since $f$ is not of
concave KPP class. Our starting point is the following 
observation.
\begin{lemma}\label{lemma1} Under assumptions of Theorem \ref{thm1} we have
  \begin{equation}
    \label{eq:gradient-ineq}
    V(t)=\frac{v_0^2}{\kappa}\int\limits_\Omega f(T(x,y,t))\frac{dxdy}{H}
\ge \kappa\int\limits_\Omega |\nabla T(x,y,t)|^2\frac{dxdy}{H}.
  \end{equation}
\end{lemma}
{\bf Proof.}  The equality in (\ref{eq:gradient-ineq}) is obtained
simply by integrating (\ref{eq:1.1}) over $\Omega$ using the boundary
conditions (\ref{eq:0.2}) and (\ref{eq:0.3}) or (\ref{eq:0.4}), and
mean-zero condition  (\ref{eq:shear-mean-zero}) on advection.  To get
the inequality we multiply (\ref{eq:1.1}) by $T$ and integrate over
$\Omega$ to get
\[
\int\limits_\Omega TT_t\frac{dxy}{H}+
\kappa\int\limits_\Omega|\nabla T(x,y,t)|^2\frac{dxdy}{H}=\frac{v_0^2}{\kappa}
\int\limits_\Omega Tf(T)\frac{dxdy}{H}\le V(t).
\]
This implies (\ref{eq:gradient-ineq}) since $T_t\ge 0$ and $0\le T\le 1$.
$\Box$

As a warm-up, we now prove a simple and general proposition, which
already provides a glimpse of some of the ideas which we will use to
obtain more sophisticated results.  Namely, we show that for any
divergence-free velocity $u(x,y)$ satisfying mild regularity
conditions (it does {\it not} have to be a shear flow), and solution $T$
satisfying (\ref{eq:1.2}), the burning rate is bounded below by
$Cv_0.$ 
\begin{prop}
\label{general}
Let T(x,y,t) be a solution of (\ref{eq:0.1}) with the boundary
conditions (\ref{eq:0.2}) and either (\ref{eq:0.3}) or (\ref{eq:0.4}).
Assume that $u(x,y) \in C^1([0,H]\times \Rm)$ satisfies
(\ref{eq:mean-zero}) and (\ref{eq:incomp}), and that non-linearity
$f(T)$ satisfies (\ref{eq:nonlin}).  Let the initial data $T_0(x,y)$
satisfy (\ref{eq:1.2.1}), (\ref{eq:in-decay}) and (\ref{eq:1.3a}).
Then there exists a constant $C,$ depending only on the parameters
$\zeta$ and $f_0,$ such that
\[ V(t) \geq C(\zeta, f_0) v_0 \]
with the constants $\zeta$ and $f_0$ defined in (\ref{eq:cond1}).
\end{prop}
{\bf Proof.} The proof is similar to the proof of Lemma 2 in \cite{CKOR}.
We can find $y$ such that 
\[ 
\int\limits_{\Rm} |\nabla T(x,y)|^2dx\leq\frac{3}{H}\int\limits_0^H 
\int\limits_{\Rm} |\nabla T(x,y')|^2 dxdy', 
\]
and
\[ 
\int\limits_{\Rm} f(T(x,y))dx\leq\frac{3}{H} 
\int\limits_0^H \int\limits_{\Rm} f(T(x,y'))dxdy'. 
\]
Then we can find $x_1 \leq x_2$ such that $T(x_1, y) =
\theta_4+\zeta,$ $T(x_2, y) = \theta_1 -\zeta,$ $T(x,y) \in
[\theta_1-\zeta, \theta_4+\zeta]$ if $x_1 \leq x \leq x_2$ (see
(\ref{eq:cond1}) for the definition of $\theta_{1,4}$). Then
\[ 
\int\limits_{\Rm} |\nabla T(x,y)|^2dx\geq 
\frac{\theta_4 -\theta_1+2\zeta}{|x_2-x_1|} 
\] 
and
\[ 
\int\limits_{\Rm}f(T(x,y))dx \geq f_0 |x_2-x_1|. 
\]
Therefore we have
\[ 
\sqrt{ \int\limits_{\Rm}\int\limits_0^H |\nabla T|^2 \, dxdy
  \int\limits_{\Rm} \int\limits_0^H f(T)\, dxdy} \geq C f_0^{1/2} \zeta H.
\]
Hence we obtain
\[ 
\int\limits_{\Rm} \int\limits_0^H \left[ \kappa |\nabla T|^2 +
  \frac{v_0^2}{\kappa} f(T) \right]\,\frac{dxdy}{H} \geq C f_0^{1/2}
\zeta v_0.
 \] 
Then Proposition \ref{general} follows from Lemma \ref{lemma1}.$\Box$

We now return to the shear flows.  To obtain more precise bounds
involving advection velocity $u(y),$ we will bound from below
in terms of $u(y)$ either the integral of $f(T)$ or the $L^2$-norm of
$|\nabla T|$, and use Lemma \ref{lemma1}. The general plan in
\cite{CKOR} was to integrate over all axis in $x$, obtaining an
equation with an explicit term $u(y)$ in it.
We were able to bound the rest of the terms from above by a combination
of $\int f(T)$ and $\int |\nabla T|^2$ after averaging in $y$ and $t$
to control $T_t$ and $\Delta T$. 

An additional twist we need here is
to reduce our consideration to the region in space where the reaction
actually takes place. In the case of ignition non-linearity, there is
no reaction for sufficiently low temperatures. Similarly, for the
Arrhenius type non-linearity, reaction is very weak near $T=0$. On the
technical side, restriction of consideration to some finite time
dependent domain $D$ with $T$ in appropriate range will mandate
additional averaging in $x$ to control all terms by $\int f(T)$ and
$\int |\nabla T|^2.$ We will identify a region $D$ in $x$ such that on
one hand the temperature has a certain drop over this region and on
the other for every $x\in D$ there is some $y\in I_j$ such that
reaction is bounded away from zero at the point $(x,y)$. This will
provide us with two alternatives for each $x\in D$: either reaction
is uniformly bounded away from zero for that $x$ or temperature drops
by a certain amount on the interval $x\times I_j$. In the first
case $\int_{I_j}f(T)dy$ will have to be large and
in the second $\int_{I_j}|T_y|^2dy$ will be
bounded from below.  Then we will integrate (\ref{eq:1.1}) over $x\in
D$ at a fixed time $t$.  That will relate $u(y)$ to some terms
involving $V(t)$ and $\Delta T$.  We will additionally average
both in $x$ and $y$, which will bring $\Delta T$ into a form that can
be bound by a combination of integrals of $f(T)$ and $|\nabla T|^2$.
That will be possible since these have to be large on $D$ as explained
above.  Finally we will use Lemma \ref{lemma1} to finish the proof.

In order to define the region $D$ where much of the reaction takes place
let us fix $\theta_4>\theta_3>\theta_2>\theta_1,$ where $\theta_4,$ $\theta_1$
are as in (\ref{eq:cond1}).
Let $I_j$ be an interval on which (\ref{eq:def-Ij}) holds with $u(y)>0$
(the case of $I_j$ where $u(y) <0$ is similar). 
We fix time $t>0$ and choose two points $x_0$ and $x_1$:
\begin{eqnarray}
  &&x_0=\inf\left\{x:\hbox{ for any $x'>x$ there exists $y\in I_j$ 
s.t. $T(x',y,t)\le \theta_4$}\right\}\nonumber\\
&&x_1=\sup\left\{x:\hbox{ $x>x_0$ and for any $x'\in(x_0,x)$ 
there exists $y\in I_j$ s.t. $T(x',y,t)\ge \theta_1$}\right\}.\nonumber
\end{eqnarray}
In other words, for any $x\in[x_0,x_1]=D$ there exists $y\in I_j$ such
that $T(x,y,t)\in[\theta_1,\theta_4]$, and hence $f(T(x,y))\ge f_0$.
Note that $x_0$ is well-defined and finite since $T(x,y,t)\to 0,1$ as
$x\to\pm\infty$ uniformly in $x$ because of (\ref{eq:later-decay}). The
definition of $x_0$ implies that 
\begin{equation}\label{eq:Tatx_0}
T(x_0,y)\ge\theta_4~~\hbox{ for all $y$}
\end{equation}
and thus $x_1$ is well-defined. Moreover, 
\begin{equation}\label{eq:Tatx_1}
\hbox{$T(x_1,y)\le\theta_1$ for all
$y$. }
\end{equation}

In preparation for multiple averaging in $y$ that will be performed to
control $\Delta T$ let us introduce the function $G(h,\xi)$:
 \begin{eqnarray*}
G(h,\xi)=\left\{\begin{matrix}
{\displaystyle{1-\frac{|\xi|}{h}}, & |\xi| \leq h\cr
0, & |\xi|>h \cr}
\end{matrix}\right.
\end{eqnarray*}
that corresponds to the following averaging in $y$:
\begin{equation}\label{aver}
\frac{1}{h_j} \int\limits_0^{h_j}
d \delta \int\limits_{c_j-\delta}^{c_j+\delta} p(y) \,dy =
\int\limits_{c_j-h_j}^{c_j+h_j} G(h_j, y-c_j)p(y) \,dy
\end{equation}
for a test function $p(y)$. The two integrations when applied
to $\Delta T$ are required to get rid of derivatives of $T$.
Observe that the function $G(h,\xi)$ has the following properties
\begin{equation}
 \label{eq:3.11}
0\le G(h,\xi)\le 1, \,\,\,\,
 G(h, \xi) \geq \frac{1}{2} ~~\hbox{for}~~
\xi\in[-\frac{h}{2},\frac{h}{2}]. 
\end{equation}
We note that (\ref{eq:Tatx_0}) and (\ref{eq:Tatx_1}) imply that
at the two ends of the interval $[x_0,x_1]$ we have
\begin{eqnarray}
&&  \int\limits_{c_j-h_j}^{c_j+h_j}dy G(h_j,y-c_j)u(y)T(x_0,y)\ge \theta_4F_j
\nonumber\\
&&\int\limits_{c_j-h_j}^{c_j+h_j}dy G(h_j,y-c_j)u(y)T(x_1,y)\le \theta_1F_j
\nonumber
\end{eqnarray}
with
\begin{eqnarray}
  F_j=\int\limits_{c_j-h_j}^{c_j+h_j}dy G(h_j,y-c_j)u(y).\nonumber
\end{eqnarray}
In preparation for averaging in $x,$ 
 we choose $\eta_0$ and $\eta_1$ so that
\begin{eqnarray}
  \label{eq:eta12}
&&\eta_0=\inf_{\xi}\left\{\xi>0:~\int\limits_{c_j-h_j}^{c_j+h_j}dy
G(h_j,y-c_j)u(y)T(x_0+\xi,y)= \theta_3F_j\right\}\nonumber\\
&&\eta_1=\inf_{\xi}\left\{\xi>0:~\int\limits_{c_j-h_j}^{c_j+h_j}dy
G(h_j,y-c_j)u(y)T(x_1-\xi,y)= \theta_2F_j\right\}.\nonumber
\end{eqnarray}
We remark that $x_0<x_0+\eta_0<x_1-\eta_1<x_1$.
Now we are ready to average (\ref{eq:1.1}). Given $\alpha\in(0,\eta_0)$
and $\beta \in(0,\eta_1)$ we integrate (\ref{eq:1.1}) in
$x\in(x_0+\alpha,x_1-\beta)$ and in $y$ according to (\ref{aver}):
\begin{eqnarray}
  \label{eq:ineq1}
  &&\int\limits_{x_0+\alpha}^{x_1-\beta}dx\int\limits_{c_j-h_j}^{c_j+h_j}dy
G(h_j,y-c_j)T_t-\kappa\int\limits_{x_0+\alpha}^{x_1-\beta}dx
\int\limits_{c_j-h_j}^{c_j+h_j}dy
G(h_j,y-c_j)T_{yy}\\&&
+\kappa\int\limits_{c_j-h_j}^{c_j+h_j}dy
G(h_j,y-c_j)[T_x(x_0+\alpha,y)-T_x(x_1-\beta,y)]\nonumber\\
&&\ge
\int\limits_{c_j-h_j}^{c_j+h_j}dy
G(h_j,y-c_j)u(y)[T(x_0+\alpha,y)-T(x_1-\beta,y)].\nonumber
\end{eqnarray}
We dropped the integral of $f(T)$ on the right side which resulted in
the inequality in (\ref{eq:eta12}). The reason that our averagings in
$x$ and $y$ are different is that while the width $h_j$ is a prescribed
quantity, we have no a priori control over $\eta_0$ and $\eta_1$.
Therefore our bounds may not involve them, and we employ different
estimates when averaging in $x$.  First we estimate the integral of $T_{yy}$ on the
left side of (\ref{eq:eta12}).
\begin{lemma}\label{lemma2} There exists a universal constant $C>0$ such that the
 following estimate holds for every 
  $x\in[x_0,x_1]$
\begin{equation}\label{eq:lemma2} 
\left|\kappa\int\limits_{c_j-h_j}^{c_j+h_j}dy
G(h_j,y-c_j)T_{yy}(x,y)\right|\le 
Cf_0^{-1/2}\zeta^{-1}\frac{\kappa }{v_0 h_j}\left[\frac{v_0^2}{\kappa}
\int\limits_{c_j-h_j}^{c_j+h_j}dyf(T(x,y))+
\kappa\int\limits_{c_j-h_j}^{c_j+h_j}dy|\nabla T|^2(x,y)\right].
\end{equation}
\end{lemma}
{\bf Proof.}  We use (\ref{aver}) to rewrite the left side of
(\ref{eq:lemma2}) for a fixed $x\in[x_0,x_1]$ as
\begin{equation}\label{eq:usedelta}
\left|~\int\limits_{c_j-h_j}^{c_j+h_j}dy
G(h_j,y-c_j)T_{yy}\right|=\frac{1}{h_j}\left| T(c_j+h_j) -2 T(c_j) + T(c_j-h_j)\right|\le \frac{2}{h_j}\delta_j[T](x)  
\end{equation}
with $\displaystyle\delta_j[T](x)=\sup_{y\in I_j}T(x,y)-\inf_{y\in I_j}T(x,y)$.
Note that because of our choice of $x_0$ and $x_1$, given any
$x\in(x_0,x_1),$ we may find $y'$ such that
$T(x,y')\in[\theta_1,\theta_4]$. Then we may find $y_1$, $y_2$ such that
for any $y\in[y_1,y_2]$ we have $T(x,y)\in(\theta_1-\zeta,\theta_4+\zeta)$,
and, moreover,
\begin{equation}
  \label{eq:perepad}
  |T(x,y_2)-T(x,y_1)|=\min(\zeta,\delta_j[T](x)).
\end{equation}
Then (\ref{eq:cond1}) implies that
\[
  \int\limits_{c_j-h_j}^{c_j+h_j}f(T(x,y))dy\ge f_0|y_2-y_1|.
\]
Applying the Cauchy-Schwartz inequality we also obtain
\[
\int\limits_{c_j-h_j}^{c_j+h_j}|\nabla T(x,y)|^2dy
\ge\frac{|T(x,y_2)-T(x,y_1)|^2}{|y_2-y_1|}.
\]
Multiplying these two inequalities
we obtain
\begin{eqnarray*}
&&\sqrt{\int\limits_{c_j-h_j}^{c_j+h_j}f(T(x,y))dy
\int\limits_{c_j-h_j}^{c_j+h_j}|\nabla T(x,y)|^2dy}
\ge \sqrt{f_0}|T(x,y_2)-T(x,y_1)|\\&&\ge \frac{f_0^{1/2} \zeta }{4}
\left| T(c_j+h_j)-2T(c_j) + T(c_j-h_j)\right|
\end{eqnarray*}
because of (\ref{eq:usedelta}) and (\ref{eq:perepad}), and since
$\zeta<1/2$. Then (\ref{eq:lemma2}) follows.$\Box$

Furthermore,   because of our choice of $\eta_0$, $\eta_1$, we have for any $\alpha\in[0,\eta_0]$ and
$\beta\in[0,\eta_1]$
\begin{eqnarray}
  \label{eq:ubd}
  \int\limits_{c_j-h_j}^{c_j+h_j}dy
G(h_j,y-c_j)u(y)[T(x_0+\alpha,y)-T(x_1-\beta,y)]\ge (\theta_3-\theta_2)F_j.
\end{eqnarray}
We use (\ref{eq:lemma2}), (\ref{eq:ubd}) and positivity of $T_t$ to
rewrite (\ref{eq:ineq1}) as
\begin{eqnarray}
  \label{eq:ineq2}
 &&
C\int\limits_{x_0+\alpha}^{x_1-\beta} dx\int\limits_{c_j-h_j}^{c_j+h_j}dy
\left[T_t+\frac{\kappa}{v_0 h_j}\left[\kappa|\nabla T|^2+
\frac{v_0^2}{\kappa}f(T(x,y))\right]\right]\\&&+
\kappa\int\limits_{c_j-h_j}^{c_j+h_j}dy
G(h_j,y-c_j)[T_x(x_0+\alpha,y)-T_x(x_1-\beta,y)]\ge
(\theta_3-\theta_2) F_j.\nonumber
\end{eqnarray}
In order to deal with the integral term on the second line that
involves $T_x$ we average (\ref{eq:ineq2}) in $\alpha\in(0,\eta_0)$ and
$\beta\in (0,\eta_1)$:
\[
\frac{1}{\eta_0\eta_1}\int\limits_0^{\eta_0}d\alpha\int\limits_0^{\eta_1}d\beta
\cdot
\]
to get
\begin{eqnarray}
  \label{eq:ineq3}
 &&
C
\int\limits_{x_0}^{x_1} dx\int\limits_{c_j-h_j}^{c_j+h_j}dy
\left[T_t+\frac{\kappa }{v_0 h_j}\left[\kappa|\nabla T|^2+
\frac{v_0^2}{\kappa}f(T(x,y))\right]\right]
\\&&+
\frac{\kappa}{\eta_0}\int\limits_{c_j-h_j}^{c_j+h_j}dy
G(h_j,y-c_j)[T(x_0+\eta_0,y)-T(x_0,y)]\nonumber\\&&+
\frac{\kappa}{\eta_1}\int\limits_{c_j-h_j}^{c_j+h_j}dy
G(h_j,y-c_j)[T(x_1-\eta_1,y)-T(x_1,y)]\ge (\theta_3-\theta_2)F_j.\nonumber
\end{eqnarray}
We bound now the term involving $\eta_0$ in (\ref{eq:ineq3}) as follows.
\begin{lemma}\label{lemma2.1}
There exists a universal constant $C>0$ such that
\begin{eqnarray}
\label{eq:lemma2.1}
 \left|\frac{\kappa}{\eta_0}\int\limits_{c_j-h_j}^{c_j+h_j}dy
G(h_j,y-c_j)[T(x_0+\eta_0,y)-T(x_0,y)]\right|\le
C \kappa (\theta_4-\theta_3)^{-2}  
\int\limits_{c_j-h_j}^{c_j+h_j}dy\int\limits_{x_0}^{x_0+\eta_0}dx
|\nabla T|^2.
\end{eqnarray}
\end{lemma}
{\bf Proof.}
The proof of this estimate is based on two observations. First, we have
\[ \left|\frac{1}{\eta_0}\int\limits_{c_j-h_j}^{c_j+h_j}dy 
G(h_j,y-c_j)[T(x_0+\eta_0,y)-T(x_0,y)]\right|\le
\frac{2 h_j}{\eta_0} \]
because of (\ref{eq:3.11}).
Second, we have 
\begin{equation}
  \label{eq:bd2}
  \int\limits_{c_j-h_j}^{c_j+h_j}dy\int\limits_{x_0}^{x_0+\eta_0}dx|\nabla T|^2\ge
C(\theta_4-\theta_3)^2 \frac{h_j}{\eta_0}.
\end{equation}
This bound is established as follows.
Recall that because of our choice of $\eta_0$ we have
\begin{eqnarray}
  \label{eq:flux1}
\int\limits_{c_j-h_j}^{c_j+h_j}dy
u(y)G(h_j,y-c_j)[T(x_0,y)-T(x_0+\eta_0,y)]\ge (\theta_4-\theta_3)F_j .
\end{eqnarray}
Furthermore, recall that $\displaystyle \frac{\|u\|_{\infty,j}}{2}\le |u(y)|
\le \|u\|_{\infty,j}$ on the interval $(c_j-h_j,c_j,+h_j)$ and thus
(\ref{eq:flux1}) implies that 
\begin{eqnarray}
  \label{eq:flux2}
  \int\limits_{c_j-h_j}^{c_j+h_j}dy
G(h_j,y-c_j)|T(x_0,y)-T(x_0+\eta_0,y)|\ge
\frac{(\theta_4-\theta_3)F_j}{\|u\|_{\infty,j}}\ge C (\theta_4-\theta_3 ) h_j.\nonumber
\end{eqnarray}
Then we obtain using the Cauchy-Schwartz inequality:
\begin{eqnarray}
  \label{eq:flux3}
&&  C(\theta_4-\theta_3) h_j\le\left(\,
 \int\limits_{c_j-h_j}^{c_j+h_j}dy
G^2(h_j,y-c_j)\right)^{1/2}
\left(\,
 \int\limits_{c_j-h_j}^{c_j+h_j}dy|T(x_0,y)-T(x_0+\eta_0,y)|^2\right)^{1/2}
\nonumber\\
&&\le Ch_j^{1/2}\left(\,
 \int\limits_{c_j-h_j}^{c_j+h_j}dy\left|
\int\limits_{x_0}^{x_0+\eta_0}dxT_x(x,y)\right|^2\right)^{1/2}\le Ch_j^{1/2}
\left(\eta_0\int\limits_{c_j-h_j}^{c_j+h_j}dy\int\limits_{x_0}^{x_0+\eta_0}dx
T_x^2(x,y)\right)^{1/2}\nonumber
\end{eqnarray}
and (\ref{eq:bd2}) follows. $\Box$

A bound similar to (\ref{eq:lemma2.1}) 
holds for the integral involving $\eta_1$ in 
(\ref{eq:ineq3}).  We use these two estimates in (\ref{eq:ineq3}) to get
\begin{eqnarray}
  \label{eq:ineq4}
&&\int\limits_{x_0}^{x_1} dx\int\limits_{c_j-h_j}^{c_j+h_j}dy
\left[T_t+\frac{\kappa }{v_0 h_j}\left[\kappa|\nabla T|^2+
\frac{v_0^2}{\kappa}f(T(x,y))\right]\right]+
{\kappa}
\int\limits_{c_j-h_j}^{c_j+h_j}dy\int\limits_{x_0}^{x_1}dx
|\nabla T|^2\\
&&\ge CF_j
\ge C'\int\limits_{c_j-h_j}^{c_j+h_j}dy|u(y)|.
\nonumber 
\end{eqnarray}
A similar estimate holds also for the intervals $I_j$, on which
$u(y)<0$.  The only difference would be that at the first step of
obtaining the analog of (\ref{eq:ineq1}) one has to drop the integral
involving $T_t$ and not that of $f(T)$. The rest of the estimates
still hold.  We use Lemma \ref{lemma1} in (\ref{eq:ineq4}) to get
after summation over all intervals $I_j$:
\[ V(t)\ge C\sum_j\left(1+\frac{l}{h_j}\right)^{-1}
\int\limits_{c_j-h_j}^{c_j+h_j}|u(y)|\frac{dy}{H}, \]
with $l=\kappa/v_0$. Finally, we can always add $v_0$ to the right
hand side by Proposition~\ref{general}.  This finishes the proof of
Theorem \ref{thm1}.

\section{Percolating flows}\label{sec:percol}

\begin{figure}
  \centerline{
  \psfig{figure=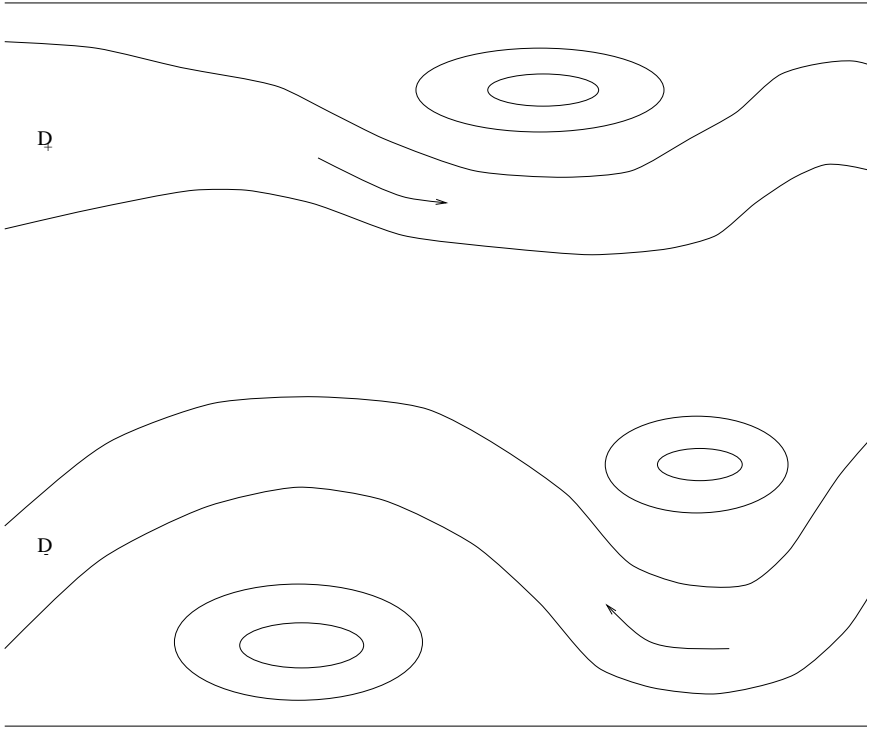,width=14cm} 
     }
  \caption{Streamlines of $u(x,y)$.}
  \label{fig1}
\end{figure}

We now consider equation
\begin{eqnarray}
  \label{eq:perc-main}
  &&T_t+u(x,y)\cdot\nabla T=\kappa\Delta T+\frac{v_0^2}{\kappa}f(T)\\
&&T(x,y,0)=T_0(x,y) \nonumber
\end{eqnarray}
in a more general class of flows, which we call ``percolating''.  By
this we mean that there exist tubes of streamlines of the advecting
velocity $u(x,y)$, which connect $x=-\infty$ and $x=+\infty$ in either
direction, as depicted on Figure \ref{fig1}. 
We assume that the flow has zero mean (\ref{eq:mean-zero})
and hence such tubes of 
streamlines will go in both directions.  
More precisely, let us assume that there
exist regions $D_j^+$ and $D_j^-$, $j=1,\dots N$ such that each of
them is bounded by the streamlines of $u(x,y)$, and the projection of
each streamline of $u(x,y)$, contained in either $D_j^+$ or $D_j^-$,
onto the $x$-axis covers the whole real line (these projections need
not be one-to-one, however). We denote by $D_\pm$ the union of all
$D_j^\pm$ respectively. 

We will further assume that the velocity $u(x,y)$ is periodic in
space. Then it is known \cite{Ber-Hamel,jxin-1,jxin-2} that for
ignition nonlinearity (\ref{eq:ignit}) there exist periodic traveling
fronts. They have the form $T(x-ct,x,y)$ and are periodic in the last
two variables and monotonically decreasing in the first one.  These
solutions satisfy our main condition
\begin{equation}
  \label{eq:per-posit}
  \pdr{T}{t}(x,y,t)\ge 0.
\end{equation}
Our results may be generalized in a straightforward manner to
non-periodic percolating flows as long as initial data satisfies
(\ref{eq:per-posit}):
\begin{equation}\label{eq:perc-init}
\kappa\Delta T_0+\frac{v_0^2}{\kappa}f(T_0)-u\cdot\nabla T_0\ge 0.
\end{equation}
However, we restrict our attention to periodic $u(x,y)$ to simplify
the presentation.

\begin{figure}
  \centerline{
  \psfig{figure=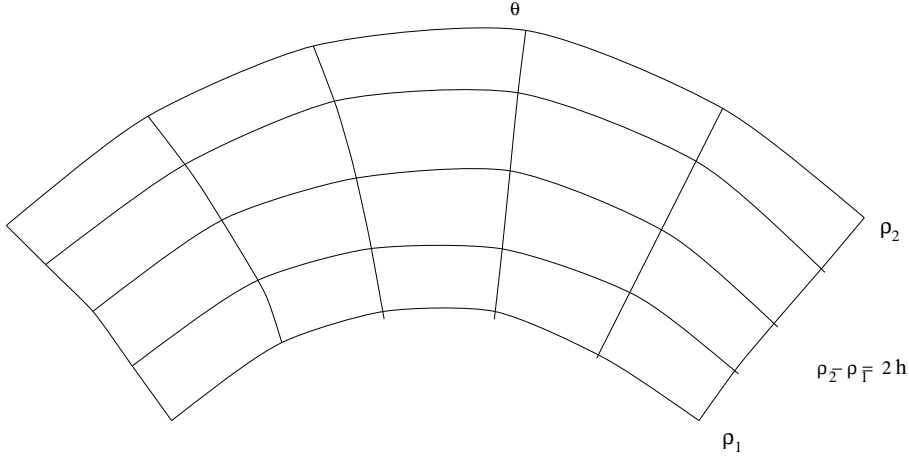,width=12cm} 
     }
  \caption{Curvilinear coordinates $(\rho,\xi)$.}
  \label{fig-coord}
\end{figure}

We assume that the streamlines in $D_j^\pm$ are sufficiently regular,
so that inside each $D_j^\pm$ there exists a one-to-one $C^2$ change
of coordinates $(x,y) \rightarrow (\rho, \xi),$ such that $\rho$ is
constant on the streamlines, while $\xi$ is an orthogonal coordinate
for $\rho$ (with a slight abuse of notation we shall use the same
notation $(\rho,\xi)$ in all $D_j^\pm,$ although these coordinates may
not be defined globally).  Moreover, $u\cdot\nabla\xi>0$ in $D_j^+$,
while $u\cdot\nabla\xi<0$ in each $D_j^-$. The variable $\rho$ varies
in the interval $[c_j^\pm-h_j^\pm, c_j^\pm+h_j^\pm]$, while $\xi$
varies in $(-\infty, \infty)$ in the set $D_j^{\pm}$. See Figure
\ref{fig-coord} for a sketch of coordinates $(\rho, \xi).$ The square
of the length element inside each set $D_j^{\pm}$ is given by
\[
dx^2+dy^2=E_1^2(\rho,\xi)d\rho^2+E_2^2(\rho,\xi)d\xi^2.
\]
 
We assume that the functions $E_{1,2}$ are bounded from below:
\begin{equation}
  \label{eq:5.0}
  C^{-1}\le E_{1,2}(\rho,\xi)
\end{equation}
uniformly on all $D_j^\pm.$ Moreover, the function 
\begin{equation}\label{eq:5.00}
\omega(\rho,\xi)=\frac{E_2(\rho,\xi)}{E_1(\rho,\xi)}
\end{equation}
satisfies the following bounds:
\begin{equation}
  \label{eq:5.1}
  C^{-1}\le|\omega(\rho,\xi)|\le C,~~
\left|\pdr{\omega}{\xi}(\rho,\xi)\right|,
\left|\frac{\partial\omega}{\partial\rho}(\rho,\xi)\right|
\le\frac{C}{h_j^{\pm}} \,\,\,\hbox{on}\,\,\,D_j^\pm \,\,\,\hbox{respectively},
\end{equation}
with $2h_j^\pm$ being the absolute value of the difference of the
values of $\rho$ on the two components of the boundary $\partial
D_j^{\pm}$ (recall that $D_j^{\pm}$ are bounded by two streamlines of
$u(x,y)$). Finally we assume that the flux density $E_1u$ does not
oscillate too much on the set $D_j^\pm$:
\begin{equation}
  \label{eq:perc-oscill}
  \frac{\|E_1u\|_{\infty,j}}{2}\le E_1(\rho,\theta)|u(\rho,\theta)|\le
\|E_1u\|_{\infty,j}~\hbox{ for }~(x,y)\in D_j^\pm,~\|E_1u\|_{\infty,j}
=\sup_{(\rho,\xi)\in D_j^\pm}|E_1(\rho,\xi)u(\rho,\xi)|.
\end{equation}
Note that $E_1(\rho,\xi)|u(\rho,\xi)|$ is independent of $\xi$. In particular 
if 
\[
u(x,y)=UH\nabla^\perp\Psi=UH\left(\pdr{\Psi}{y},-\pdr{\Psi}{x}\right)
\]
with $|\nabla\Psi|\le C/H$ we may choose $\rho=H\Psi(x,y)$ so that
$E_1=\frac{1}{H|\nabla\Psi|}$. Then we have $E_1|u|=U$, so that
(\ref{eq:perc-oscill}) holds automatically and (\ref{eq:5.0}) also
holds for $E_1$. Other conditions on the streamlines may be also
easily restated in terms of the stream function $\Psi(x,y)$.

We do not make any assumptions on the behavior of the streamlines of
$u(x,y)$ outside the regions $D_+$ and $D_-$. In particular, there may
be pockets of still fluid, streamlines may be closed, etc. (see Figure
\ref{fig1}).

Then we have the following Theorem.
\begin{thm}\label{thm2} Let $T(x,y,t)$ be a solution of (\ref{eq:perc-main})
  with the boundary conditions (\ref{eq:0.2}) and (\ref{eq:0.3}) or
  (\ref{eq:0.4}), with the initial data $T_0(x,y)$ satisfying
  (\ref{eq:perc-init}), and nonlinearity $f(T)$ satisfying
  (\ref{eq:nonlin}).  Let each of the sets $D_j^{\pm}$ be of the form
  $D_j^{\pm}=\left\{\rho\in[c_j-h_j,c_j+h_j]\right\}$. Then under the
  assumptions (\ref{eq:5.0}) and (\ref{eq:5.1})
on the streamlines of the flow $u(x,y)\in C^1(\Omega)$, we have
\begin{eqnarray}
  \label{eq:5.2}
V(t)\ge C\left(v_0+\sum_{D_j^\pm}
\left(1+\frac{l}{h_j}\right)^{-1}
\int\limits_{c_j-{h_j}}^{c_j+{h_j}}|u(\rho, \xi)|
E_1(\rho, \xi)
\frac{d\rho}{H}\right)
\end{eqnarray}
for all $t>0$. Here $l=\kappa/v_0$
and the constant $C$ in (\ref{eq:5.2}) depends only
on the function $f(T)$ and the constants appearing in (\ref{eq:5.0})
and (\ref{eq:5.1}).
\end{thm}
Note that the integrals on the right side of (\ref{eq:5.2}) give
fluxes of $u(x,y)$ through the tubes of the streamlines.
As in the shear case the pre-factor $(1+l/h_j)^{-1}$ agrees with the
homogenization limit \cite{CKOR,HPS}. 

Recall that traveling fronts for periodic flows have the form 
\begin{equation}
\label{eq:per-travel}
T(x,y,t)=U(x-ct,x,y)
\end{equation}
with the function $U(s,x,y)$ being periodic in the last two variables.
It was shown in \cite{Ber-Hamel,jxin-1,jxin-2} that in the ignition
nonlinearity case (\ref{eq:ignit}) such traveling fronts exist and
their speed $c_*$ is unique. Their existence for the KPP
nonlinearities (\ref{eq:kpp}) was shown recently in \cite{Ber-Hamel}
with $c\ge c_u$, $c_u$ being the minimal traveling front speed. 
The following analog of Corollary~\ref{cor1} holds for percolating
flows, which we formulate separately for the convenience of the
reader.
\begin{corollary}\label{cor2}
  Let $f(T)$ be either of the ignition nonlinearity type
  (\ref{eq:ignit}), or of the KPP type (\ref{eq:kpp}). Let also
  $U(x-ct,x,y)$ be a traveling wave-type solution of (\ref{eq:perc-main}),
  periodic in the second two variables. Then there exists a constant
  $C>0$ that depends only on the function $f$ and on the constants
  appearing in (\ref{eq:5.0}) and (\ref{eq:5.1}) such that
 \[ c\ge C\left( v_0+\sum_j\left(1+\frac{l}{h_j}\right)^{-1}
\int\limits_{c_j-h_j}^{c_j+{h_j}}|u(\rho, \xi)|
E_1(\rho, \xi)
\frac{d\rho}{H}\right). \]
\end{corollary}
{\bf Proof.} 
Corollary \ref{cor2} follows from Theorem \ref{thm2} as follows.  Let
$U(x-ct,x,y)$ be a periodic traveling front solution of
(\ref{eq:perc-main}) such that 
\[
U(s,x+L,y)=U(s,x,y)
\]
and let $\tau=L/c$. Then we have
\begin{eqnarray*}
&&\frac{1}{\tau}\int
\limits_0^{\tau}V(t)dt=\frac{1}{\tau}\int\limits_0^{\tau}dt\int\limits_{\Omega}
(-cU_s(x-ct,x,y))\frac{dxdy}{H}=
\frac{1}{\tau}\int\limits_0^{L}dt'\int\limits_{\Omega}[-U_s(x-t',x,y)]
\frac{dxdy}{H}\\
&&=
-\frac{c}{L}\int\limits_0^{L}dt'\int\limits_{\Omega}\left[\frac{d}{dx}
U(x-t',x,y)
-U_x(x-t',x,y)\right]\frac{dxdy}{H}\\&&=c+
\frac{c}{L}\int\limits_0^{L}dt'\int\limits_{\Omega}U_x(x-t',x,y)\frac{dxdy}{H}
=c+\frac{c}{L}\int\limits_{\Omega}\frac{dxdy}{H}[U(x,x+L,y)-U(x,x,y)]=c.~\Box
\end{eqnarray*}

{\it Remark.}  We note that Corollary \ref{cor2} implies a lower bound on the
effective diffusivity \cite{BLP} in the homogenization regime. 
Recall that solutions of the advection-diffusion equation
(\ref{eq:perc-main}) with $f(T)=0$, and with advection of the form
$\displaystyle u(x,y)=\frac{U}{\varepsilon}
v\left(\frac{x}{\varepsilon},\frac{y}{\varepsilon}\right)$ 
with $v({x,y})$ periodic, converge as $\varepsilon\to 0$ to
the solution $\bar T$ of the homogenized problem
\[
\pdr{\bar T}{t}=\kappa_{ij}^*\frac{\partial^2 \bar T}{\partial x_i\partial
  x_j},~~\bar T(x,y,0)=T_0(x,y),~(x_1,x_2)=(x,y).
\]
The effective diffusivity $\kappa^*$ is a complicated functional of
advection $Uv({\bf y})$. Explicit bounds on $\kappa^*$ in terms of the
magnitude $U$ of advection are easy to obtain in the shear case, when
$\kappa^*$ may be found explicitly, and effective diffusivity in the
direction of the flow $\kappa_{xx}^*\sim U^2$.  Using results of
\cite{HPS} one may deduce from Corollary \ref{cor2} that this bound
applies also to periodic percolating flows, despite the fact that
no explicit expression for $\kappa^*$ is known in this case.

Corollary \ref{cor2} and the stability results of \cite{jxin-2}
imply the analog of Theorem \ref{thm0.1}.
\begin{thm}\label{cor3}
  Let the initial data  $T_0(x,y)$ for the equation (\ref{eq:perc-main})
satisfy the decay
  to $0,1$ conditions (\ref{eq:in-decay}) and (\ref{eq:1.3a}). 
 Then both for KPP and ignition
non-linearities we have
 \[  \lim_{t\to\infty}\langle V\rangle_t\ge 
 C\left( v_0+\sum_j\left(1+\frac{l}{h_j}\right)^{-1}
\int\limits_{c_j-{h_j}}^{c_j+{h_j}}|u(\rho, \xi)|
E_1(\rho, \xi)
\frac{d\rho}{H}\right), \]
where constant $C$ depends only on the non-linearity $f$ and on the
constants appearing in (\ref{eq:5.0}) and (\ref{eq:5.1}).
\end{thm}

Now we turn to the proof of Theorem \ref{thm2}.  The proof is a
modification of the proof of Theorem \ref{thm1}. We will again utilize
Lemma \ref{lemma1} as well as averaging along the streamlines of
$u(x,y)$ in order to bound the arising averages of $\Delta T$ in terms
of integrals of $f(T)$ and $|\nabla T|^2$.  The additional technical
difficulties are due to the fact that two natural geometries of the
problem - streamlines for the advective term and Euclidean coordinates
for Laplacian - are in harmony in the case of shear flows, but at odds
in the case of more general percolating flows. Moreover, while in the
case of shear flows we gave the proof that we felt was simplest, here
we will use the approach which is slightly more involved; however, it
is better adapted for the application to cellular flows in the
following sections.

Let us consider a region $D_j^+=\left\{(\rho,\xi):
  \rho\in[c_j-h_j,c_j+h_j]\right\}$ with $u\cdot\nabla\xi>0$.
Introduce notation $k_j(\rho) = G(h_j, \rho-c_j)E_1(\rho, \xi)|u(\rho,
\xi)|$ ($k_j$ does not depend on $\xi$ by incompressibility of
$u(x,y)$) and, similarly to the shear case,
\[ 
F_j = \int\limits_{c_j -h_j}^{c_j+h_j} k_j(\rho)\, d \rho. 
\]
 We can find $\xi_1>\xi_0$ such that 
\[ \int\limits_{c_j-h_j}^{c_j+h_j} k_j(\rho) T(\rho, \xi_0) \, 
d \rho = \theta_4 F_j, \,\,\,
\int\limits_{c_j-h_j}^{c_j+h_j} k_j(\rho) T(\rho, \xi_1) \, 
d \rho = \theta_1 F_j,  \]
and for every $\xi \in [\xi_0, \xi_1]$ we have 
\[ \theta_1 F_j \leq \int\limits_{c_j-h_j}^{c_j+h_j} k_j(\rho)
 T(\rho, \xi) \, d \rho \leq \theta_4 F_j. \]
Let us denote by $\tilde{D}_j^+$ the region bounded by the curves
$\rho = c_j \pm h_j$ and $\xi= \xi_0, \xi_1.$ This region depends
on time, but we will suppress this dependence in notation.
Theorem \ref{thm2} will follow from the following
\begin{thm}
\label{advest0} 
 Let $T(x,y,t)$ be a solution of (\ref{eq:perc-main})
  with the boundary conditions (\ref{eq:0.2}) and (\ref{eq:0.3}) or
  (\ref{eq:0.4}), with the initial data $T_0(x,y)$ satisfying
  (\ref{eq:perc-init}), and nonlinearity $f(T)$ satisfying
  (\ref{eq:nonlin}). Then under the
  assumptions (\ref{eq:5.0}) and (\ref{eq:5.1})
on the streamlines of the flow $u(x,y)\in C^1(\Omega)$, we have
for every time $t$
\[ \int\limits_{\tilde{D}_j^+}\left(T_t+\kappa|\nabla T|^2 + 
\frac{v_0^2}{\kappa} f(T) \right) \, dxdy \geq 
C(\zeta, f_0) \left(1+ \frac{l}{h}\right)^{-1} (\theta_4-\theta_1)^3 F_j  \]
with the constants $\zeta$, $f_0$ and $\theta_{1,4}$ defined by
(\ref{eq:cond1}).
\end{thm}
\it Remark. \rm A similar lower bound on the integral over $D^+_j$ is easier
to obtain, and also suffices to prove Theorem~\ref{thm2}.
 We chose, however, to formulate Theorem~\ref{advest0} in this stronger 
version since this is what we will need when dealing with 
cellular flows. \\

\noindent {\bf Proof.}
Let $\theta_3 = \theta_4 - (\theta_4-\theta_1)/3,$ and $\theta_2 = \theta_1+(\theta_4-\theta_1)/3.$ 
Then we may choose $\eta_0$ and
$\eta_1$ similarly to the shear case, namely
\begin{eqnarray*}
&& \eta_0=\inf_{\eta}\left\{\eta>0:~\int\limits_{c_j-h_j}^{c_j+h_j}d\rho
k_j(\rho) T(\rho,\xi_0+\eta)= 
\theta_3F_j\right\}\nonumber\\
&& \eta_1=\inf_{\eta}\left\{\eta>0:~\int\limits_{c_j-h_j}^{c_j+h_j}d\rho
k_j(\rho)T(\rho,\xi_1-\eta)= 
\theta_2 F_j\right\}
\end{eqnarray*}
 We choose
$\alpha\in(0,\eta_0)$ and $\beta\in(0,\eta_1),$ integrate
(\ref{eq:perc-main}) over $\xi\in (\eta_0+\alpha,\eta_1-\beta)$ and average
in $\rho$ with the kernel $G$ to get 
\begin{eqnarray}
  \label{eq:perc-est1}
 && \int\limits_{\xi_0+\alpha}^{\xi_1-\beta}d\xi\int\limits_{c_j-h_j}^{c_j+h_j}d\rho
G(h_j,\rho-c_j)E_1E_2T_t-\kappa
\int\limits_{\xi_0+\alpha}^{\xi_1-\beta}d\xi\int\limits_{c_j-h_j}^{c_j+h_j}d\rho
G(h_j,\rho-c_j)E_1E_2[T_{xx}+T_{yy}]\\
&&  \ge
\int\limits_{c_j-h_j}^{c_j +h_j} 
k_j(\rho) [T(\rho,\xi_0+\alpha)-
T(\rho, \xi_1-\beta)] d \rho \nonumber \\ && \ge \frac{1}{3}(\theta_4-\theta_1)F_j.
\nonumber
\end{eqnarray}
 In (\ref{eq:perc-est1}), we dropped the term involving $f(T)$ on the right-hand side, as we did in
the shear case.  We now
look at the term involving the Laplacian:
\begin{eqnarray*}
&& \int\limits_{\xi_0+\alpha}^{\xi_1-\beta}d\xi
\int\limits_{c_j-h_j}^{c_j+h_j}d\rho
G(h_j,\rho-c_j)E_1(\rho,\xi)E_2(\rho,\xi)[T_{xx}(\rho,\xi)+T_{yy}(\rho,\xi)]\\
&&=
\frac{1}{h_j}\int\limits_0^{h_j}d\delta
\int\limits_{\xi_0+\alpha}^{\xi_1-\beta}d\xi
\int\limits_{c_j-\delta}^{c_j+\delta}
d\rho E_1(\rho,\xi)E_2(\rho,\xi)[T_{xx}(\rho,\xi)+T_{yy}(\rho,\xi)].
\end{eqnarray*}
Let us denote $D_{\delta \alpha \beta}$ the region bounded by 
coordinate curves $\rho = c_j \pm \delta$ and $\xi= \xi_0+\alpha,$ 
$\xi = \xi_1 -\beta.$ 
Using Green's formula,
we rewrite the last two integrations for fixed $\delta$ as
\begin{eqnarray}
  &&\int\limits_{\xi_0+\alpha}^{\xi_1-\beta}d\xi
\int\limits_{c_j-\delta}^{c_j+\delta}
d\rho E_1(\rho,\xi)E_2(\rho,\xi)[T_{xx}(\rho,\xi)
+T_{yy}(\rho,\xi)]=\nonumber\\
&& = \int \int_{D_{\delta \alpha \beta}} \Delta T \, dxdy =
\int_{\partial D_{\delta \alpha \beta}} \frac{\partial T}{\partial n}ds 
\nonumber \\
&&=\int\limits_{\xi_0+\alpha}^{\xi_1-\beta}d\xi
\left[\omega(c_j+\delta,\xi)\pdr{T}{\rho}(c_j+\delta,\xi)-
\omega(c_j-\delta,\xi)\pdr{T}{\rho}(c_j-\delta,\xi)\right]\nonumber
\\&&
+\int\limits_{c_j-\delta}^{c_j+\delta}
d\rho\left[{\omega^{-1}}(\rho,\xi_1-\beta)\pdr{T}{\xi}(\rho,\xi_1-\beta)-
{\omega^{-1}}(\rho,\xi_0+\alpha)\pdr{T}{\xi}(\rho,\xi_0+\alpha)\right].
\label{eq:perc-lapl}
\end{eqnarray}
 We used the definition (\ref{eq:5.00}) of the function
$\omega(\rho,\xi)=E_2(\rho,\xi)/E_1(\rho, \xi)$ in the last step.
The average of the term on the first line in (\ref{eq:perc-lapl}) is
estimated by the following Lemma, an analog of Lemma~\ref{lemma2}

\begin{lemma}\label{lemma3} There exists a constant $C>0$ such that we have
for all $\xi \in [\xi_0,\xi_1]$
  \begin{eqnarray}
    \label{eq:lem3}
&&\frac{\kappa}{h_j}\left|
\int\limits_0^{h_j}d\delta
\left[\omega(c_j+\delta,\xi)\pdr{T}{\rho}(c_j+\delta,\xi)-
\omega(c_j-\delta,\xi)\pdr{T}{\rho}(c_j-\delta,\xi)\right]\right|\\&&\le
Cf_0^{-1/2}\zeta^{-1}\frac{\kappa }{v_0 h_j}\left[
\frac{v_0^2}{\kappa}\int\limits_{c_j-h_j}^{c_j+h}d\rho
 E_1E_2f(T) +\kappa
\int\limits_{c_j-h_j}^{c_j+h_j}d\rho
 E_1E_2|\nabla T|^2\right],\nonumber
\end{eqnarray}
where the constant $C$ depends only on constants in the bounds (\ref{eq:5.1}) and (\ref{eq:5.0}). 
\end{lemma}
{\bf Proof.} 
We will show that
\begin{eqnarray}
  \label{eq:perc-est-fixed-xi}
 \left|
\int\limits_0^{h_j}d\delta
\omega(c_j+\delta,\xi)\pdr{T}{\rho}(c_j+\delta,\xi)\right|\le 
 \frac{Cf_0^{-1/2}\zeta^{-1}}{v_0}\left[
\kappa \int\limits_{c_j}^{c_j+h_j}d\rho T_\rho^2(\rho,\xi)+
\frac{v_0^2}{\kappa}\int\limits_{c_j}^{c_j+h_j}d\rho f(T(\rho,\xi))\right]
\end{eqnarray}
for all $\xi\in[\xi_0,\xi_1]$. 
A similar estimate holds for the second term on
the left side of (\ref{eq:lem3}). 
By definition of $\xi_0$ and $\xi_1$ for every
$\xi\in[\xi_0,\xi_1]$ there exists $\rho_0\in (c_j-h_j,c_j+h_j)$ such that
$T(\rho_0,\xi)\in[\theta_1,\theta_4]$. Then given $\xi\in[\xi_0,\xi_1]$ we
have two possibilities.  First, assume that
$T\in(\theta_1-\zeta,\theta_4+\zeta)$ for all
$\rho\in(c_j-h_j,c_j+h_j)$. Then we have for such $\xi$
\[
\int\limits_{c_j}^{c_j+h_j}d\rho f(T(\rho,\xi))
\ge Cf_0\int\limits_{c_j}^{c_j+h_j}d\rho
\omega^2(\rho,\xi),
\]
which implies that
\begin{eqnarray*}
&&\left|
\int\limits_0^{h_j}d\delta
\omega(c_j+\delta,\xi)\pdr{T}{\rho}(c_j+\delta,\xi)\right|\le 
\left(\int\limits_{c_j}^{c_j+h_j}d\rho\omega^2(\rho,\xi)\right)^{1/2}
\left(\int\limits_{c_j}^{c_j+h_j}d\rho T_\rho^2(\rho,\xi)\right)^{1/2}\\
&&
\le
Cf_0^{-1/2}\left(\int\limits_{c_j}^{c_j+h_j}d\rho f(T(\rho,\xi))\right)^{1/2}
\left(\int\limits_{c_j}^{c_j+h_j}d\rho T_\rho^2(\rho,\xi)\right)^{1/2}\\
&&\le \frac{Cf_0^{-1/2}}{v_0}\left[
\kappa \int\limits_{c_j}^{c_j+h_j}d\rho T_\rho^2(\rho,\xi)+
\frac{v_0^2}{\kappa}\int\limits_{c_j}^{c_j+h_j}d\rho f(T(\rho,\xi))\right].
\end{eqnarray*}
The other case for a given $\xi$
is that the temperature $T(\rho,\xi)$ drops out of the range
$(\theta_1-\zeta,\theta_4+\zeta)$ for some $\rho$. Then we may find
$\rho'\in(c_j-h_j,c_j+h_j)$ such that $T(\rho',\xi)=\theta_1-\zeta$ or
$T(\rho',\xi)=\theta_4+\zeta$, and, moreover,
$T\in(\theta_1-\zeta,\theta_4+\zeta)$ for all $\rho$ between $\rho'$
and $\rho_0$. Then we have $|T(\rho',\xi)-T(\rho_0,\xi)|\ge\zeta$ and
hence
\begin{eqnarray*}
\left(\int\limits_{c_j}^{c_j+h_j}d\rho f(T(\rho,\xi))\right)^{1/2}
\left(\int\limits_{c_j}^{c_j+h_j}d\rho T_\rho^2(\rho,\xi)\right)^{1/2}
\ge \left(\frac{\zeta^2}{|\rho'-\rho_0|}f_0|\rho'-\rho_0|\right)^{1/2}
\ge \zeta f_0^{1/2}.
\end{eqnarray*}
Then the estimate (\ref{eq:perc-est-fixed-xi}) also holds in that case since
\begin{eqnarray*}
  \left|
\int\limits_0^{h_j}d\delta
\omega(c_j+\delta,\xi)\pdr{T}{\rho}(c_j+\delta,\xi)\right|\le
C
\end{eqnarray*}
for all $\xi$, as can be seen from integrating by parts in $\rho$ and using
(\ref{eq:5.1}). Since $\zeta \leq 1/2,$  (\ref{eq:perc-est-fixed-xi}) holds in
both cases.  $\Box$

In order to bring the term in the second line of (\ref{eq:perc-lapl})
into a form convenient for analysis we average (\ref{eq:perc-lapl}) in
$\alpha$ and $\beta$. Let us consider an estimate on the second
summand after averaging in $\alpha$ over $[0,\eta_0];$ the other summand
is treated similarly.
\begin{lemma}\label{lemma5}
  There exists a constant $C,$ depending only on constants in bounds
  (\ref{eq:5.0}) and (\ref{eq:5.1}), such that
\begin{eqnarray}
  \label{eq:lem5}
  &&\left|\frac{\kappa}{\eta_0 h_j}
\int\limits_0^{h_j}d\delta
\int\limits_0^{\eta_0}d\alpha\int\limits_{c_j-\delta}^{c_j+\delta}d\rho
\omega^{-1}(\rho,\xi_0+\alpha)\pdr{T}{\xi}(\rho,\xi_0+\alpha)\right|\\
&&\le C\left[f_0^{-1/2}  
(\theta_4-\theta_1)^{-1} {\rm min}((\theta_4-\theta_1), \zeta)^{-1}
 \frac{\kappa }{v_0 h_j}\int\limits_{\tilde{D}_j^+}dxdy
\left[\frac{v_0^2}{\kappa}f(T(x,y))
+\kappa|\nabla T(x,y)|^2\right] \right.\nonumber \\&& 
\left.+(\theta_4-\theta_1)^{-2}\kappa \int\limits_{\tilde{D}_j^+}dxdy
|\nabla T(x,y)|^2\right]. \nonumber
\end{eqnarray}
\end{lemma}
{\bf Proof.}
We have to estimate the following expression
\begin{eqnarray}
  &&\frac{1}{\eta_0}
\int\limits_0^{h_j}d\delta
\int\limits_0^{\eta_0}d\alpha\int\limits_{c_j-\delta}^{c_j+\delta}d\rho
\omega^{-1}(\rho,\xi_0+\alpha)\pdr{T}{\xi}(\rho,\xi_0+\alpha)\nonumber\\
\label{eq:lem5-transf}&&=
\frac{1}{\eta_0}
\int\limits_0^{h_j}d\delta
\int\limits_{c_j-\delta}^{c_j+\delta}d\rho
[\omega^{-1}(\rho,\xi_0+\eta_0){T}(\rho,\xi_0+\eta_0)-
\omega^{-1}(\rho,\xi_0){T}(\rho,\xi_0)]\\
&&-
\frac{1}{\eta_0}
\int\limits_0^{h_j}d\delta
\int\limits_0^{\eta_0}d\alpha\int\limits_{c_j-\delta}^{c_j+\delta}d\rho
\pdr{\omega^{-1}}{\xi}(\rho,\xi_0+\alpha){T}(\rho,\xi_0+\alpha).\nonumber
\end{eqnarray}
The first term in (\ref{eq:lem5-transf}) is evidently bounded by
\begin{eqnarray}\label{eq:lem5-pf-0}
\frac{1}{\eta_0}
\left|
\int\limits_0^{h_j}d\delta
\int\limits_{c_j-\delta}^{c_j+\delta}d\rho
[\omega^{-1}(\rho,\xi_0+\eta_0){T}(\rho,\xi_0+\eta_0)-
\omega^{-1}(\rho,\xi_0){T}(\rho,\xi_0)]\right|\le\frac{Ch_j^2}{\eta_0}.
\end{eqnarray}
Furthermore, we claim that
\begin{equation}
  \label{eq:lem5-pf-1}
  \int\limits_{\tilde{D}_j^+} 
dxdy|\nabla T(x,y)|^2\ge\frac{C(\theta_4-\theta_3)^2 h_j}{\eta_0}.
\end{equation}
This is shown exactly as the estimate (\ref{eq:bd2}) in the proof of
Lemma~\ref{lemma2.1} for
the shear case, given the assumption 
(\ref{eq:5.1}) on $\omega.$ 
We combine (\ref{eq:lem5-pf-0}) and
(\ref{eq:lem5-pf-1}) to obtain
\begin{eqnarray}
&&\frac{1}{\eta_0}
\left|
\int\limits_0^{h_j}d\delta
\int\limits_{c_j-\delta}^{c_j+\delta}d\rho
[\omega^{-1}(\rho,\xi_0+\eta_0){T}(\rho,\xi_0+\eta_0)-
\omega^{-1}(\rho,\xi_0){T}(\rho,\xi_0)]\right|\nonumber\\
\label{eq:lem5-pf-7}
&&\le
C(\theta_4-\theta_3)^{-2}h_j\int\limits_{\tilde{D}_j^+}
dxdy|\nabla T(x,y)|^2.\label{eq:lem5-pf-2}
\end{eqnarray}
Next we estimate the second term in (\ref{eq:lem5-transf}):
\begin{equation}
\label{eq:lem5-pf-3}  
\frac{1}{\eta_0 }
\left|
\int\limits_0^{h_j}d\delta
\int\limits_{\xi_0}^{\xi_0+\eta_0}d\xi
\int\limits_{c_j-\delta}^{c_j+\delta}d\rho
\pdr{\omega^{-1}}{\xi}(\rho,\xi){T}(\rho,\xi)\right|
\le Ch_j
\end{equation}
because of (\ref{eq:5.1}). 
Notice that if $\eta_0 \leq h_j,$ then Lemma~\ref{lemma5}
follows directly from (\ref{eq:lem5-pf-1}). Hence we can assume  
$\eta_0>h_j.$   
The following final Lemma allows us to finish the proof of Lemma \ref{lemma5}.

\begin{lemma}\label{lemma6}
Assume that $\eta_0>h_j.$ Then there exists a constant 
$C,$ which depends only on the constants 
in bounds (\ref{eq:5.1}) and (\ref{eq:5.0}), such that 
\[ \int\limits_{\tilde{D}_j^+} f(T) 
\,dxdy \int\limits_{\tilde{D}_j^+}  |\nabla T|^2 \,dxdy 
\geq C f_0 
(\theta_4-\theta_1)^2 
\left({\rm min} ((\theta_4-\theta_1),\zeta/2)\right)^2 h_j^2. \]
\end{lemma} 
\it Remark. \rm This lemma is much easier to prove with $D_j^+$
instead of $\tilde{D}_j^+,$ basically the argument of 
Proposition~\ref{general}
applies. The simpler version is also sufficient for the proof 
of Theorem~\ref{thm2}, but we need this stronger version for 
the proof of Theorem~\ref{advest0}. \\

\noindent {\bf Proof.} Let us denote 
$|S|$ the Lebesgue measure of a measurable set $S.$ 
We know that
\[ \int\limits_{c_j-h_j}^{c_j+h_j} k_j(\rho)(T(\rho, 
\xi_0)-T(\rho, \xi_1))\, d\rho =
(\theta_4-\theta_1)F_j. \]
Consider two sets
\[ S_1 = \{ \rho \left| \right. T(\rho, \xi_0) - T(\rho, \xi_1) < 
(\theta_4-\theta_1)
/2 \}; \,\,\,S_2 = [c_j-h_j,c_j+h_j] \setminus S_1. \]
Then 
\[ \int\limits_{S_1} k_j(\rho) (T(\rho, \xi_0)-T(\rho, \xi_1))\,d \rho 
\leq (\theta_4-\theta_1)F_j/2, \]
and hence 
\[  \int\limits_{S_2} k_j(\rho) (T(\rho, \xi_0)-T(\rho, \xi_1))\,d \rho 
\geq (\theta_4-\theta_1)F_j/2. 
\] 
Since $T(\rho, \xi_0)-T(\rho, \xi_1)
\leq 1,$ and $k_j(\rho)/F_j \leq 2h_j^{-1}$ by the properties of $G$
and (\ref{eq:perc-oscill}) we must have $|S_2| \geq
(\theta_4-\theta_1)h_j/4.$ Choose $\rho_0 \in S_2$ such that
\[ 
\int\limits_{\xi_0}^{\xi_1} d\xi f(T(\rho_0,\xi)) \leq 
\frac{3}{|S_2|}\int\limits_{\tilde{D}^+_j} f(T)\,dxdy 
\]
and 
\[  
\int\limits_{\xi_0}^{\xi_1} d\xi T_\xi^2(\rho_0,\xi) \leq 
\frac{3}{|S_2|}\int\limits_{\tilde{D}^+_j} T_\xi^2 \,d\rho d\xi \leq
\frac{C}{|S_2|}\int\limits_{\tilde{D}^+_j} |\nabla T|^2\,dxdy, 
\]
with $C$ depending only on the constants in (\ref{eq:5.1}).  Since
$\rho_0 \in S_2,$ we have $T(\rho_0, \xi_0)-T(\rho_0, \xi_1)\geq
(\theta_4-\theta_1)/2.$ Assume first that there exists $\xi \in [\xi_0,
\xi_1]$ where $T(\rho_0, \xi) \in [\theta_1-\zeta/2,
\theta_4+\zeta/2].$ Then we can find $\xi_0\leq \xi_2 <\xi_3\leq
\xi_1$ such that for every $\xi \in [\xi_2, \xi_3],$ $T(\rho_0, \xi)
\in [\theta_1-\zeta, \theta_4+\zeta],$ and
\[ 
T(\rho_0, \xi_2)-T(\rho_0, \xi_3) \geq {\rm
  min}((\theta_4-\theta_1)/2, \zeta/2).
\]
Therefore, 
\[ 
\int\limits_{\xi_0}^{\xi_1} f(T) d \xi \geq f_0 (\xi_3-\xi_2), 
\]
while 
\[ 
\int\limits_{\xi_0}^{\xi_1} T_\xi^2  d \xi \geq  
\frac{({\rm min}((\theta_4-\theta_1)/2, \zeta/2))^2}{\xi_3-\xi_2}. 
\]
Combining these estimates and the definition of $\rho_0,$ we find
\[ 
\int\limits_{\tilde{D}_j^+} f(T)\,dxdy 
\int\limits_{\tilde{D}^+_j} |\nabla T|^2\,dxdy \geq
Cf_0 ({\rm min}((\theta_4-\theta_1)/2, \zeta/2))^2 |S_2|^2, 
\]
which proves the lemma in this case. 

The other case we have to consider is that for every $\xi \in [\xi_0,
\xi_1],$ $T(\rho_0, \xi) \notin [\theta_1-\zeta/2, \theta_4+\zeta/2]$.
Then we will be able to find a point $\xi\in[\xi_0,\xi_1]$ so that 
there is a drop in temperature along the curve $\xi=\hbox{const}$.
Assume $T(\rho_0, \xi)>\theta_4+\zeta/2$ for all $\xi,$ the other case
being similar. By mean value theorem, we also have that for every $\xi
\in [\xi_0, \xi_1]$ there exists $\rho_1$ such that $T(\rho_1, \xi)
\in [\theta_1, \theta_4],$ and so $T(\rho_1, \xi)\leq \theta_4.$
Similarly to the above, we can find $\tilde{\xi}$ such that
\[ 
\int\limits_{c_j-h_j}^{c_j+h_j} d\rho f(T(\rho,\tilde{\xi})) \leq
\frac{3}{\xi_1-\xi_0}\int\limits_{\tilde{D}^+_j} f(T)\,dxdy \leq
\frac{3}{\eta_0}\int\limits_{\tilde{D}^+_j} f(T)\,dxdy
\] 
and
\[  
\int\limits_{c_j-h_j}^{c_j+h_j} d\rho T_\rho^2(\rho,\tilde{\xi}) \leq 
\frac{3}{\xi_1-\xi_0}\int\limits_{\tilde{D}_j^+} T_\rho^2 \,d\rho d\xi \leq
\frac{3}{\eta_0}\int\limits_{\tilde{D}_j^+} T_\rho^2 \,d\rho d\xi \leq
\frac{C}{\eta_0}\int\limits_{\tilde{D}_j^+} |\nabla T|^2\,dxdy. 
\]
An argument identical to the one we used in the previous case establishes that
\[ 
\int\limits_{c_j-h_j}^{c_j+h_j} f(T) d \rho \geq f_0 (\rho_1-\rho_2), 
\]
while 
\[ 
\int\limits_{c_j-h_j}^{c_j+h_j} T_\rho^2  d \rho \geq  
\frac{\zeta^2}{4(\rho_1-\rho_2)} 
\]
for some $c_j+h_j \geq \rho_1>\rho_2 \geq c_j-h_j.$ 
This implies that 
\[ 
\int\limits_{\tilde{D}_j^+} f(T)\,dxdy 
\int\limits_{\tilde{D}_j^+} |\nabla T|^2\,dxdy \geq
Cf_0\zeta^2 \eta_0^2 \geq C f_0\zeta^2 h_j^2. 
\]
Thus Lemma \ref{lemma6} is proven. $\Box$

Lemma \ref{lemma6} and the assumptions (\ref{eq:5.0}) and (\ref{eq:5.1})
imply that
\[
\int\limits_{\tilde{D}_j^+}{dxdy}
\left[\frac{v_0}{\kappa}f(T)+\frac{\kappa}{v_0}|\nabla
  T|^2\right]\ge Cf_0^{1/2}(\theta_4-\theta_1) 
{\rm min}((\theta_4-\theta_1)/2, \zeta/2) h_j.
\]
Therefore we deduce from (\ref{eq:lem5-pf-3}) that
\begin{eqnarray}
\label{eq:lem5-pf-4}  
&&\frac{1}{\eta_0 h_j}
\left|
\int\limits_0^{h_j}d\delta
\int\limits_{\xi_0}^{\xi_0+\eta_0}d\xi
\int\limits_{c_j-\delta}^{c_j+\delta}d\rho
\pdr{\omega^{-1}}{\xi}(\rho,\xi){T}(\rho,\xi)\right|
\\ 
&&\le Cf_0^{-1/2}(\theta_4-\theta_1)^{-1}
 {\rm min}((\theta_4-\theta_1)/2, \zeta/2)^{-1}
 \frac{1}{v_0 h_j}\int\limits_{\tilde{D}_j^+}{dxdy}
\left[\frac{v_0^2}{\kappa}f(T)+\kappa|\nabla
  T|^2\right].\nonumber
\end{eqnarray}
Then (\ref{eq:lem5-transf}), (\ref{eq:lem5-pf-7}) and (\ref{eq:lem5-pf-4})
imply Lemma \ref{lemma5}. $\Box$

Averaging of (\ref{eq:perc-lapl}) in $\beta\in[0,b_1]$ produces terms
similar to the left side of (\ref{eq:lem5}), which are bounded by the
same quantity.  Now we put together our estimates on $\tilde{D}_j^+$, that is,
equations (\ref{eq:perc-est1}) and
(\ref{eq:perc-lapl}), and Lemmas \ref{lemma3} and \ref{lemma5} :
\begin{eqnarray*}
  &&\int\limits_{\tilde{D}_j^+}dxdy
G(h_j,\rho-c_j)T_t+C \zeta^{-1}f_0^{-1/2}
\frac{\kappa}{v_0 h_j}\int\limits_{\tilde{D}^+_j}
dxdy
\left[\frac{v_0^2}{\kappa}f(T)+\kappa|\nabla T|^2\right]+ \\&&
 Cf_0^{-1/2}(\theta_4-\theta_1)^{-1} 
{\rm min}((\theta_4-\theta_1), \zeta)^{-1} 
\frac{\kappa}{v_0 h_j}\int\limits_{\tilde{D}_j^+}dxdy
\left[\frac{v_0^2}{\kappa}f(T)+\kappa
|\nabla T|^2\right] \\ && 
+C(\theta_4-\theta_1)^{-2}\kappa \int\limits_{\tilde{D}_j^+}dxdy|\nabla T|^2
\ge C(\theta_4-\theta_1)F_j.
\end{eqnarray*}
Therefore we obtain
\begin{eqnarray*}
  && \left(1+\frac{\kappa}{v_0h_j}\right)\int\limits_{\tilde{D}_j^+} dxdy
\left[T_t+\frac{v_0^2}{\kappa}f(T)+\kappa|\nabla T|^2\right]\ge
C(\zeta, f_0) (\theta_4-\theta_1)^3 F_j  \\ && \geq
C(\zeta, f_0) (\theta_4-\theta_1)^3\int\limits_{c_j-h_j}^{c_j+h_j}E_1(\rho,\xi)|u(\rho,\xi)|d\xi,
\end{eqnarray*}
which proves Theorem \ref{advest0}. $\Box$

The estimate on the regions $D_j^-$ with $u\cdot\nabla\xi<0$ is
similar. The only essential difference is that in inequality
(\ref{eq:perc-est1}) we drop the term involving $T_t$ but keep the one
with $f(T)$, which does not make any difference in the final result.
Then Theorem \ref{thm2} follows after summation over all $D_j^\pm$ 
from Lemma~\ref{lemma1} and Proposition~\ref{general}. $\Box$

\section{Cellular flows: the main result}\label{sec:cell-main}

Now we consider (\ref{eq:0.1}) in a cellular flow:
\begin{equation}
  \label{eq:celleq}
  T_t+u\cdot\nabla T=\kappa\Delta T+\frac{v_0^2}{\kappa}f(T),
\end{equation}
 that is, flow with closed
streamlines, and establish a lower bound for the burning rate. 
For simplicity, we limit our consideration to one typical representative class
of cellular flows, given by the stream function 
\begin{equation}\label{stream}
 \psi(x,y)= UH \sin \frac{x}{H} \sin \frac{y}{H} 
\end{equation}
on the strip $(-\infty,\infty)\times [0,\pi H]$ (for convenience from now on 
$H$ will be the width of the strip divided by $\pi$). 
The flow $u(x,y)$ is given by
\begin{equation}\label{eq:stream1}
u(x,y)=U\nabla^\perp\psi(x,y)=U\left(\sin \frac { x}H\cos\frac {y}H, 
-\cos\frac {x}H\sin\frac {y}H\right).
\end{equation}
The streamlines inside a period cell are depicted on Figure \ref{fig:cell}.
The results we prove can be extended in a direct
way to the periodic cellular flows of more general form.

We will further assume that the Peclet number is larger than one:
\begin{equation}
  \label{eq:peclet}
  \hbox{Pe}=\frac{UH}{\kappa}\ge 1
\end{equation}
and the size of the cell is larger than the laminar front width:
\begin{equation}
  \label{eq:cellsize}
  \frac{l}{H}\le 1, ~~
l=\frac{\kappa}{v_0}.
\end{equation}
Conditions (\ref{eq:peclet}) and (\ref{eq:cellsize}) are natural for
flows of large amplitude and for thin fronts. Moreover, our results
may be easily adapted to the other regimes, where (\ref{eq:peclet}) and
(\ref{eq:cellsize}) are violated.  We introduce also the turnover time
$\tau_u$, and the chemical reaction time $\tau_c$:
\begin{equation}
  \label{eq:turnover}
  \tau_u=\frac{H}{U},~~\tau_c=\frac{\kappa}{v_0^2}. 
\end{equation}
It turns out that the ratio $\tau_u/\tau_c$ is the crucial parameter 
for burning in the cellular flows. 

Finally we assume that $T(x,y,t)$ satisfies the usual boundary
conditions and that
\[
\pdr{T}{t}(x,y,t)\ge 0.
\]
As we noted previously, this condition is satisfied as long as it
holds initially.
\begin{thm}\label{cellthm} 
  Let $T(x,y,t)$ be a solution of (\ref{eq:1.1}) with the boundary
  conditions (\ref{eq:0.2}) and either (\ref{eq:0.3}) or
  (\ref{eq:0.4}), and the cellular flow given by (\ref{eq:stream1}). Let
  the initial data $T_0(x,y)$ satisfy (\ref{eq:1.2.1}),
  (\ref{eq:in-decay}) and (\ref{eq:1.3a}), and let the non-linearity
  $f(T)$ be of either ignition or general KPP type. Furthermore,
  assume that (\ref{eq:peclet}) and (\ref{eq:cellsize}) hold. Then
  we have for any time $t$
\begin{equation}
\label{perter1}
V(t)\ge\left\{
{\begin{matrix}{
\left(C_1\sqrt{\frac{\tau_c}{\tau_u}}+C_2\right)v_0, &
{\rm if}~~  
{\tau_c}\le\tau_u\cr
\left(C_1\left(\frac{\tau_c}{\tau_u}\right)^{1/5}+C_2\right)v_0, &
{\rm if}~~  
{\tau_c}\ge\tau_u\cr}\end{matrix}}\right.
\end{equation}
The constants in the inequalities depend only on the reaction $f$,
more particularly on constants $f_0,$ $\zeta,$ and
$\theta_4-\theta_1$ that appear in (\ref{eq:cond1}).
\end{thm}
\noindent \it Remark. \rm As we noted above, in order to avoid
excessive details, we chose not to formulate Theorem~\ref{cellthm} in
the exhaustive form which goes through all possible relationships
between parameters (large $\kappa$ limit, small $H$ limit, small $v_0$
limit). The reader will find it not difficult to extend the results we
prove to the above mentioned regimes.  Theorem~\ref{cellthm} is
formulated here for the range of parameters that appears to be
physically reasonable for most problems of interest.

Furthermore, we have the following corollary.
\begin{corollary}\label{cellcor} 
  Let $f(T)$ be of ignition nonlinearity type (\ref{eq:ignit}), or of
  the KPP type (\ref{eq:kpp}) and let $c$ be the speed of a
  traveling wave-type solution $T(x,y,t)=U(x-ct,x,y)$ of
  (\ref{eq:celleq}), periodic in the second two variables. Then there
  exist constants $C_{1,2}>0$ which depend only on the function $f$
  and on the constants appearing in (\ref{eq:5.0}) and (\ref{eq:5.1})
  such that
\[ c\ge\left\{
{\begin{matrix}{
\left(C_1\sqrt{\frac{\tau_c}{\tau_u}}+C_2\right)v_0, &
{\rm if}~~  
{\tau_c}\le\tau_u\cr
\left(C_1\left(\frac{\tau_c}{\tau_u}\right)^{1/5}+C_2\right)v_0, &
{\rm if}~~{\tau_c}\ge\tau_u\cr}\end{matrix}}\right. \]
\end{corollary}
Corollary \ref{cellcor} follows from Theorem \ref{cellthm} since the
traveling front profile $U(s,x,y)$ is monotonically decreasing in $s$
(see \cite{jxin-1} for the ignition case, and \cite{Ber-Hamel} for
KPP).  Then Theorem \ref{celltravelthm} follows immediately from
Corollary \ref{cellcor}, results of \cite{jxin-1} and the argument we
gave in the proof of Theorem~\ref{thm0.1}.

The proof of Theorem \ref{cellthm} is a boundary layer argument that 
proceeds, roughly, as follows. 
The temperature drops from one on the left to zero
on the right.  We will watch the temperature in the layers of width
$h$ formed by streamlines near the boundary of the cells.  The drop of
temperature in these layers may occur inside the cells or over the
diffusive interfaces. The first estimate, which we call advective,
shows how much the cell must contribute to the bulk burning rate if a
certain drop of the temperature (in the range $[\theta_1, \theta_4]$)
takes place along the streamlines inside the cell. It is reasonable to
expect that the drop over the cell will be small when advection is
strong since it mixes the fluid inside the cell quickly; in an
analytic form this intuition will translate into a large lower bound
for the burning rate if the temperature drop is significant.  The
second estimate, which we call diffusive, gives a lower bound for the
burning rate given certain drop of the temperature between the two
cells. We do expect the temperature to drop on the boundaries, and
hence the lower bound is only effective if we choose $h$ in an
appropriate way, sufficiently small.  Finally, we prove the reaction
estimate, which takes into account the total area of the region over
which the temperature drops. These estimates will be brought together to
establish the lower bound for the bulk burning rate using an appropriate 
optimization argument.

\section{Cellular flows: regularity of the streamlines}

Our first objective is to define appropriate curvilinear coordinates on the
cells, and to show that these coordinates satisfy certain technical
assumptions that we will need. The natural choice of the coordinate
$\rho$, which is constant on the streamlines, is
\begin{equation}\label{rho}
 \rho(x,y)= \frac{\psi (x,y)}{U} = H \sin \frac{x}{H} \sin \frac{y}{H}.
\end{equation}
We have certain 
freedom in the definition of the orthogonal coordinate $\xi$ along the 
streamlines : 
\begin{equation}
\label{zeta}
\nabla \xi = Q \nabla^\perp \rho,
\end{equation}
where $Q$ is some function 
which should satisfy 
\begin{equation}
\label{Qeq}
\nabla \rho\cdot \nabla Q = -Q \Delta \rho.
\end{equation} 
It is easy to compute that 
\[ \partial_\rho = x_\rho \partial_x + y_\rho \partial_y = 
\frac{\xi_y}{J}\partial_x - \frac{\xi_x}{J} \partial_y, \]
where $J = \rho_x \xi_y - \rho_y \xi_x = -Q |\nabla \rho|^2.$
Hence 
\[ 
\partial_\rho = \frac{1}{|\nabla \rho|^2}\nabla\rho\cdot\nabla,
\]
and so we have from (\ref{Qeq})
\begin{equation}\label{Qeq1}
 \frac{\partial Q}{\partial \rho} = - \frac{Q \Delta \rho}{|\nabla \rho|^2}.
\end{equation}
We will choose $Q$ so that $Q(H/2, \xi)=1;$ equation (\ref{Qeq1}) then 
allows us to define $Q$ in the region $H > \rho > 0:$
\begin{equation}
\label{QQ}
Q(\rho, \xi) = e^{-\int\limits_{H/2}^\rho \frac{\Delta \rho(h, \xi)}
{|\nabla \rho (h, \xi)|^2} \,dh}.
\end{equation}
We have the following auxiliary 
\begin{lemma}
\label{Qest}
In the region $H/2 \geq \rho >0$ we have the following bound for the
function $Q:$
\[ 0<e^{-1}\leq Q(\rho, \xi) \leq e. \]
\end{lemma}
{\bf Proof.}  Consider formula (\ref{QQ}). We have 
\[  \Delta \rho = -\frac{2\pi^2}{H^2}  \rho(x,y), \]
and 
\begin{eqnarray}
\label{gradest} 
 &&|\nabla \rho|^2 =  \left(\cos \frac{x}{H} \sin \frac{y}{H}\right)^2 +
 \left(\sin \frac{x}{H} \cos \frac{y}{H}\right)^2\nonumber \\
&&=(\sin \frac{y}{H})^2+(\sin \frac{x}{H})^2  - \frac{2\rho^2}{H^2} \geq 
2\frac{\rho}{H} (1- \frac{\rho}{H}) \geq \frac{\rho}{H}
\end{eqnarray}
in the region $\rho \leq H/2.$ 
Hence in the region  $H/2 \geq \rho >0,$ 
\[  \frac{|\Delta \rho|}
{|\nabla \rho |^2} \leq 2/H. \]
Therefore, by (\ref{QQ}), we have $e^{-1} \leq Q \leq e$ in this region. $\Box$

Recall our notation
\[ dx^2 + dy^2 = E_1^2 d \rho^2 + E_2^2 d \xi^2, \,\,\,\omega=\frac{E_2}
{E_1}. \]
The next proposition  
 summarizes some of the properties of the coordinates 
$(\rho,\xi)$ in a region of interest to us. 
\begin{prop}
\label{procoord}
For the cellular flow defined by (\ref{eq:stream1}) and coordinates
$\rho,$ $\xi$ defined by (\ref{rho}), (\ref{zeta}), and (\ref{QQ}),
the following bounds hold in the region $H/2 \geq  \rho >0:$
\begin{eqnarray}
\label{Eest}
E_{1,2} (\rho, \xi) \geq C  \\
\label{west} 0<C^{-1} \geq \omega(\rho, \xi) \leq C \\
\label{wrhoest}\left|\frac{\partial \omega}{\partial \rho}(\rho, \xi)\right| 
\leq CH^{-1} 
\\ \nonumber
\left| \frac{\partial \omega}{\partial \xi}(\rho, \xi)\right|
\leq CH^{-1}\left|\log (\rho/H)\right|.
\end{eqnarray}
\end{prop}
{\bf Proof.}  Direct computation gives that $E_1= 1/|\nabla \rho|,$
$E_2 = 1/Q|\nabla \rho|$.  Then Lemma~\ref{Qest} and the fact that
$|\nabla\rho|\le 1$ imply (\ref{Eest}).  It follows from the definition
of $\omega$ that $\omega = 1/Q$, and hence Lemma~\ref{Qest} implies
(\ref{west}).  Next,
\[ 
\left|\frac{\partial\omega}{\partial\rho}\right|=\frac{1}{Q}\frac{|\Delta\rho|}
{|\nabla \rho|^2} \leq  \frac{2e}{H}, 
\]
proving (\ref{wrhoest}). 
Finally, 
\[ 
\left|\frac{\partial\omega}{\partial\xi}\right| = 
\frac{2}{QH^2} \left| \int_{H/2}^{\rho} h \partial_\xi
\left( \frac{1}{|\nabla \rho|^2}\right) \,dh \right|. 
\]
Notice that 
\[ \partial_\xi = x_\xi \partial_x + y_\xi \partial_y = 
\frac{\rho_y}{Q|\nabla \rho|^2}\partial_x -
 \frac{\rho_x}{Q|\nabla \rho|^2} \partial_y. 
\]
A straightforward computation using (\ref{QQ}) and (\ref{gradest}) leads to 
\[ 
\left|\frac{\partial \omega}{\partial\xi}\right|\leq 
\frac{C}{H}\left| \int_{H/2}^\rho\frac{dh}{h}\right|
\leq CH^{-1}\left|\log(\rho/H)\right|. ~~~~~~\Box 
\]

\section{Cellular flows: Advective estimate}\label{sec:advect}

Let us introduce some notation. Within the cell, we will normalize
$\xi$ by letting it be zero in the negative direction of the $x$ axis
(assuming that the origin has been placed in the center of the cell).
We will denote by $L$ the value of $\xi$ in the positive direction of the
$x$ axis. In every cell, we will consider a tube of streamlines
bounded by $\rho = h$ and $\rho = 3h;$ $h$ will be always assumed to
be less than $H/6.$ We set
\[ 
k(\rho) = G(h, \rho-2h)E_1(\rho, \xi) |u(\rho, \xi)|. 
\]
The fact that $k(\rho)$ does not depend on $\xi$ is a direct corollary
of incompressibility of the flow. Moreover, with our definition of $\rho$ for
the cellular flow, we have $k(\rho)= U G(h, \rho-2h).$ We also denote
\[ 
F = \int\limits_h^{3h} k(\rho)\, d\rho. 
\]
As a corollary of Proposition~\ref{procoord}, in the strip $h \leq
\rho \leq 3h$ we have all conditions on $E_{1,2}$ and $\omega$ that
have been necessary for the percolating flow estimates. In
particular, we have
\begin{equation}
\label{omxiest}
\left|\frac{\partial\omega}{\partial\xi}\right|\leq CH^{-1}\left|
\log(\rho/H)\right|\leq Ch^{-1}, 
\end{equation}
since $h\le\rho\le H$.
We first state an estimate very similar to Theorem~\ref{advest0}. 
\begin{thm}
\label{1adv}
Assume that within one cell $C,$ there exist two values $\xi_0,$
$\xi_1$ such that for some $s_0, s_1 \in [\theta_1, \theta_4],$ $s_0 >
s_1,$ we have
\begin{eqnarray*}
\int\limits_{h}^{3h} k(\rho) T(\rho, \xi_0) \, d \rho = s_0 F, \,\,\,\, 
\int\limits_{h}^{3h} k(\rho) T(\rho, \xi_1)\, d \rho= s_1 F, \\
s_1 F \leq \int\limits_{h}^{3h} k(\rho) T(\rho, \xi) \, d \rho 
\leq s_0 F \,\,\,{\rm for}\,\,\,\xi \in [\xi_0, \xi_1].
\end{eqnarray*}
Let  $D$ be the region bounded by 
the curves $\rho=h,$ $\rho=3h,$ $\xi=\xi_0$ and $\xi=\xi_1.$ 
Then 
\begin{equation}
\label{advest1}
\int\limits_D \left(T_t+
\kappa |\nabla T|^2 + \frac{v_0^2}{\kappa} f(T) \right) \,dxdy
\geq C(\zeta,f_0)\left(1+\frac{l}{h}\right)^{-1}(s_0-s_1)^3 F.
\end{equation}
The constant $C$ in (\ref{advest1}) depends only on parameters $\zeta$
and $f_0$ of the reaction $f$, that appear in (\ref{eq:cond1}), and on
the constants in the bounds of Proposition~\ref{procoord}.
\end{thm}
{\bf Proof.} The proof is exactly the same as for
Theorem~\ref{advest0} for percolating flows.  We need only to replace
$\theta_1$ with $s_1,$ $\theta_4$ with $s_0,$ and set $c_j =2h,$
$h_j=h.$ $\Box$

The estimate (\ref{advest1}) works well if there is a significant
change of the temperature along the streamlines within the cell. But
it has a serious flaw if the temperature drops gradually and there is
little change of temperature inside any cell. The factor $(s_0-s_1)^3$
on the right hand side of the estimate makes it rather inefficient.
Our next goal is to derive an estimate which has linear
dependence on the temperature drop when the drop is small.  Our main
measurement tool for the temperature within the cell will be the
following average
\[ 
\langle T \rangle_{\xi_0} = \frac{1}{2AF} 
\int\limits_{\xi_0-A}^{\xi_0+A} d \xi 
\int\limits_h^{3h} k(\rho) T(\rho, \xi)\,d\rho. 
\]
We will normally take $A$ so that the region of averaging in $\xi$
covers about half of the width of the strip when $\xi_0=0$ or
$\xi_0=L,$ hence $A \approx H.$ Our measure of the temperature change
along the cell will be the difference between such averages for
different $\xi_0.$
\begin{thm}
\label{advest2}
Assume that for all $\xi$ in a given cell we have 
\[  
\frac{1}{F}\int\limits_h^{3h} k(\rho)T(\rho, \xi) d \rho 
\in [\theta_1, \theta_4].
\]
Then for any $\xi_0,$ $\xi_1$ we have for $A\ge h$
\begin{equation}
\label{adv2}
\int\limits_{\xi_0-A}^{\xi_1+A} \int\limits_h^{3h} \left( 
T_t+\kappa |\nabla T|^2 +
\frac{v_0^2}{\kappa}f(T)\right)E_1E_2 d\xi d\rho \geq
C(f_0, \zeta) \left( 1+\frac{l}{h}\right)^{-1} 
\left|  \langle T \rangle_{\xi_0} -\langle T \rangle_{\xi_1}\right| F, 
\end{equation}
where $C$ depends only on $f_0,$ $\zeta$ and constants in the bounds
of Proposition~\ref{procoord}.
\end{thm}
{\bf Proof.} We will consider the case where $\langle T
\rangle_{\xi_0} >\langle T \rangle_{\xi_1},$ the other case being
similar.  Let us integrate
\[ 
T_t + u\cdot \nabla T - \kappa \Delta T = \frac{v_0^2}{\kappa} f(T) 
\]
in $x$ and $y$ over the region where $\xi$ varies from $\xi_2$ to
$\xi_3,$ where $\xi_2 \in [\xi_0-A, \xi_0+A],$ $\xi_3 \in [\xi_1-A,
\xi_1+A],$ and $\rho$ takes values between $h$ and $3h$, with kernel
$G(h, \rho-2h).$ We obtain
\begin{eqnarray} 
\label{startest}
&&\int\limits_{\xi_2}^{\xi_3}d\xi \int\limits_h^{3h} 
d\rho G(h, \rho-2h)E_1E_2 T_t(\rho,\xi) - \kappa 
\int\limits_{\xi_2}^{\xi_3}d\xi \int\limits_h^{3h} 
d\rho G(h, \rho-2h)E_1E_2 \Delta T  \nonumber \\
&&\ge\int\limits_h^{3h} d \rho k(\rho) (T(\rho, \xi_2)-T(\rho, \xi_3)).
\end{eqnarray}
As usual, our goal is to control the Laplacian term. Rewrite it in a
familiar form
\begin{equation}
\label{Lapterm}
\frac{\kappa}{h} \int\limits_0^{h} d \delta \left\{ 
\int\limits_{\xi_2}^{\xi_3} d\xi \left[
\frac{\partial T}{\partial \rho}\omega(\xi, 2h+\delta)  - 
\frac{\partial T}{\partial \rho}\omega(\xi, 2h-\delta)\right] + 
\int\limits_{2h-\delta}^{2h+\delta} d\rho 
\left[\frac{\partial T}{\partial \xi}\omega^{-1}(\xi_3, \rho)
-\frac{\partial T}{\partial \xi}\omega^{-1}(\xi_2, \rho)\right] \right\}. 
\end{equation}
We carry out two more averagings in (\ref{startest}) 
\[  
\frac{1}{4A^2}  \int\limits_{\xi_0-A}^{\xi_0+A} d \xi_2
 \int\limits_{\xi_1-A}^{\xi_1+A} d \xi_3 
\]
to be able to estimate the second part of the Laplacian term
(\ref{Lapterm}). Let us start by estimating the first expression in square
brackets, more precisely,
\[  \frac{\kappa}{h}  \int\limits_0^{h} d \delta 
\int\limits_{\xi_2}^{\xi_3} d\xi  \frac{\partial T}{\partial \rho}
\omega(\xi, 2h+\delta) \]
(the other part is estimated similarly). Notice that by the 
assumption of the theorem 
for  every $\xi$ there exists $\rho$ (depending on $\xi$) such that 
$T(\rho, \xi) \in [\theta_1, \theta_4].$  
The estimate of this term now follows step by step the estimate
of the same term in the proof of 
Theorem~\ref{advest0}. We get that for every $\xi \in [\xi_0-A, \xi_1+A],$ 
\begin{equation}
\label{Lap1}
  \frac{\kappa}{h} \int\limits_0^{h} d \delta 
 \frac{\partial T}{\partial \rho}\omega(\xi, 2h+\delta) 
\leq C f_0^{-1/2}\zeta^{-1}\frac{\kappa}{v_0 h} 
\int\limits_h^{3h} \left[ \kappa |\nabla T|^2 + \frac{v_0^2}{\kappa} 
f(T)\right] \,E_1E_2 d\rho. 
\end{equation}
 The integrations in $\xi_{2,3}$  simply
average out,  they are not needed for this term.

Let us now consider the estimate of the second expression in square
brackets in (\ref{Lapterm}), more particularly the first summand (the
second one is estimated in the same way). Averaging in $\xi_2$ simply
disappears since there is no dependence on this variable, and we are
left with
\[ \frac{\kappa}{2hA} \int\limits_0^{h} d \delta
  \int\limits_{2h-\delta}^{2h+\delta} d\rho  \int\limits_{\xi_1-A}^{\xi_1+A} 
d \xi_3 \frac{\partial T}{\partial \xi}\omega^{-1}(\xi_3, \rho). \]
The following lemma is the crucial step in the proof.
\begin{lemma}
\label{lemma10}
Under conditions of Theorem~\ref{advest2} we have
\begin{equation}
\label{lem10}
\frac{\kappa}{2hA} \left| \int\limits_0^h d \delta
  \int\limits_{2h-\delta}^{2h+\delta} d\rho  \int\limits_{\xi_1-A}^{\xi_1+A} 
d \xi_3 \frac{\partial T}{\partial \xi}\omega^{-1}(\xi_3, \rho) \right| \leq 
C f_0^{-1/2} \zeta^{-1} \frac{\kappa}{h v_0} \int\limits_D 
\left[\kappa |\nabla T|^2+\frac{v_0^2}{\kappa} f(T)\right]\,dxdy,
\end{equation}
where the constant $C$ depends only on the bounds in
Proposition~\ref{procoord}.
\end{lemma}
{\bf Proof.} 
Integrating by parts and using (\ref{omxiest}), we find that
\begin{equation}\label{eq:genb0}
\left|\,\int\limits_{\xi_1-A}^{\xi_1+A} 
d\xi_3\frac{\partial T}{\partial \xi}\omega^{-1}(\xi_3, \rho) \right| \leq
C\left(1+\frac{A}{h}\right)
\end{equation}
and therefore
\begin{equation}
\label{genb}
\frac{\kappa}{2hA} \left| \int\limits_0^{h} d \delta
  \int\limits_{2h-\delta}^{2h+\delta} d\rho  \int\limits_{\xi_1-A}^{\xi_1+A} 
d \xi_3 \frac{\partial T}{\partial \xi}\omega^{-1}(\xi_3, \rho) \right| \leq
C\kappa \left( \frac{h}{A}+1 \right)\le C\kappa. 
\end{equation}
We have to consider several cases. \\
{\it Option 1.} There exists $\rho_0$ such that for all $\xi \in
[\xi_1-A, \xi_1+A],$ we have $T(\rho_0, \xi) \notin [\theta_1-\zeta/2,
\theta_4 +\zeta/2].$ By mean value theorem, for every $\xi,$ there
also exists $\rho_1(\xi)$ such that $T(\rho_1, \xi) \in [\theta_1,
\theta_4].$ Then for every $\xi$ we can find $\rho_3>\rho_2$ such that
$|T(\rho_2, \xi)-T(\rho_3, \xi)| = \zeta/2,$ and for every $\rho \in
[\rho_2, \rho_3],$ $T(\rho, \xi) \in [\theta_1-\zeta/2,
\theta_4+\zeta/2].$ Therefore, as we have seen before, for every 
$\xi \in [\xi_1-A, \xi_1+A]$ we have
\[ 
\left( \int\limits_h^{3h} d \rho T_\rho^2 \right)^{1/2} 
\left( \int\limits_h^{3h} d \rho f(T) \right)^{1/2}
\geq \frac{\zeta f_0^{1/2}}{2}. 
\]
Hence, using Cauchy-Schwartz and integrating in $\xi,$ we get 
\[ 
\int\limits_h^{3h} d \rho \int\limits_{\xi_1-A}^{\xi_1+A}d \xi E_1E_2 
\left( \kappa |\nabla T|^2 + \frac{v_0^2}{\kappa}
f(T) \right) \geq \frac{v_0 \zeta f_0^{1/2} A}{2}. 
\]
From this inequality and (\ref{genb}) our lemma follows since $A\ge h$. \\
{\it Option 2.} For every $\rho,$ there exists $\xi_4$ such that
$T(\rho, \xi_4) \in [\theta_1-\zeta/2, \theta_4+\zeta/2].$ Here we
have to consider two distinct sets of $\rho.$ First, assume that
$T(\rho, \xi) \in [\theta_1-\zeta, \theta_4+\zeta]$ for all $\xi.$
Denote the set of all such $\rho$ by $S_1.$ For $\rho \in S_1$ we have
\[ 
\int\limits_{\xi_1-A}^{\xi_1+A} f(T) d\xi 
\geq C f_0 \int\limits_{\xi_1-A}^{\xi_1+A} \omega^{-2}(\rho, \xi) d\xi,
\]
and so 
\begin{eqnarray}
\label{estimate1}
&& \left| \, \int\limits_{\xi_1-A}^{\xi_1+A}  
\frac{\partial T}{\partial \xi}\omega^{-1}(\rho, \xi) \, d\xi \right|
\leq Cf_0^{-1/2} \left(\,
\int\limits_{\xi_1-A}^{\xi_1+A} f(T) d\xi\right)^{1/2} 
\left(\, 
\int\limits_{\xi_1-A}^{\xi_1+A} T_{\xi}^2 d\xi \right)^{1/2} \\ && \leq 
Cf_0^{-1/2}v_0^{-1} \int\limits_{\xi_1-A}^{\xi_1+A}
 \left[ \kappa |\nabla T|^2 + \frac{v_0^2}{\kappa} f(T) \right]
E_1E_2 \,d\xi. \nonumber 
\end{eqnarray}
The second case is that there exists $\xi_4 \in [\xi_1-A, \xi_1+A]$
such that $T(\rho, \xi_4) \notin [\theta_1-\zeta, \theta_4+\zeta].$
Denote the set of all such $\rho$ by $S_2.$ For $\rho \in S_2,$ we can
find $\xi_5$ and $\xi_6$ such that $|T(\rho, \xi_5)-T(\rho, \xi_6)| =
\zeta/2,$ and $T(\rho, \zeta) \in [\theta_1-\zeta, \theta_4+\zeta]$
for all $\xi$ between $\xi_5$ and $\xi_6.$ In this case, similarly to
the above reasoning, we have
\[ \int\limits_{\xi_1-A}^{\xi_1+A} f(T(\rho, \xi)) d\xi
 \int\limits_{\xi_1-A}^{\xi_1+A}  T^2_\xi(\rho, \xi) d\xi
\geq C\zeta^2 f_0, \]
and so
\[ f_0^{-1/2}\zeta^{-1} v_0^{-1} 
\int\limits_{\xi_1-A}^{\xi_1+A} \left[ \kappa |\nabla T|^2 +
 \frac{v_0^2}{\kappa} f(T) \right]
E_1E_2 \,d\xi \geq C. \]
Then (\ref{eq:genb0}) implies that for $\rho\in S_2$ we have
\begin{equation}
  \label{eq:genb5}
 \left|\,\int\limits_{\xi_1-A}^{\xi_1+A} 
d\xi_3\frac{\partial T}{\partial \xi}\omega^{-1}(\xi_3, \rho) \right| \leq
C\left(1+\frac{A}{h}\right)f_0^{-1/2}\zeta^{-1} v_0^{-1} 
\int\limits_{\xi_1-A}^{\xi_1+A} \left[ \kappa |\nabla T|^2 +
 \frac{v_0^2}{\kappa} f(T) \right]
E_1E_2 \,d\xi.
\end{equation}
Combining the two estimates (\ref{estimate1}) and (\ref{eq:genb5}),
and taking into account that $A \geq h,$ we obtain the result of 
Lemma \ref{lemma10}. $\Box$
 
Now we can finish the proof of Theorem~\ref{advest2}. Taking into account 
(\ref{Lap1}) and (\ref{lem10}) we see that 
\begin{eqnarray*}
&&  \int\limits_{\xi_0-A}^{\xi_1+A}d\xi \int\limits_h^{3h} 
d\rho G(h, \rho-2h)E_1E_2 T_t(\rho,\xi) +
 Cf_0^{-1/2}\zeta^{-1} \frac{\kappa}{v_0 h}
\int\limits_{\xi_0-A}^{\xi_1+A}d\xi \int\limits_h^{3h} 
d\rho E_1E_2\left[ \kappa |\nabla T|^2 +\frac{v_0^2}{\kappa}
f(T) \right] \geq   \\ &&
\frac{1}{4A^2}\int\limits_{\xi_0-A}^{\xi_0+A} d\xi_2 
\int\limits_{\xi_1-A}^{\xi_1+A} d\xi_3
\int\limits_h^{3h} d \rho k(\rho) (T(\rho, \xi_2)-T(\rho, \xi_3)), 
\end{eqnarray*}
and this implies 
\[ 
\int\limits_{\xi_0-A}^{\xi_1+A} \int\limits_h^{3h} \left( 
T_t+\kappa |\nabla T|^2 +
\frac{v_0^2}{\kappa}f(T)\right)E_1E_2 d\xi d\rho \geq
C(f_0, \zeta) \left( 1+\frac{l}{h}\right)^{-1} 
\left|  \langle T \rangle_{\xi_0} -\langle T \rangle_{\xi_1}\right| F.
 \,\,\,\,\,\,\, \Box 
\]

\section{Cellular flows: Diffusive estimate}\label{sec:diff}

\begin{figure}
  \centerline{
  \psfig{figure=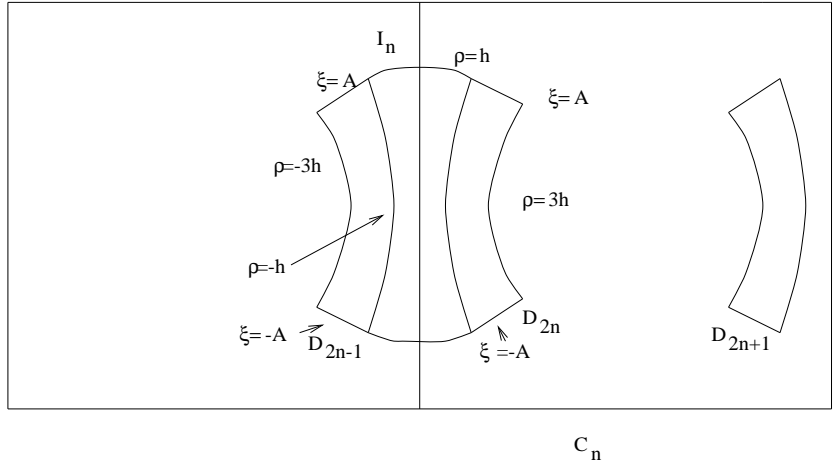,width=11cm} 
     }
  \caption{The regions $D_n$.}
  \label{fig:neighbor}
\end{figure}

Our  goal in this section is to estimate the burning rate
from below in terms of the jump of the temperature across the interface
separating two cells. We will only consider this estimate in the context
of the particular cellular flow we are studying, though it can be 
easily extended to a more general situation. Consider two neighboring
cells, which we denote $C_{n-1}$ and $C_{n}.$ We will look at two
regions $D_{2n-1} \subset C_{n-1}$ and $D_{2n} \subset C_{n}$
which are symmetric under reflection 
with respect to the line separating the cells (see Figure \ref{fig:neighbor}). 
The region $D_{2n}$ is bounded by the curves $\rho = h,$ $\rho = 3h,$
$\xi = -A,$ and $\xi = A.$ For simplicity, we choose $\rho$ and $\xi$
coordinates in cell $C_{n-1}$ so that $\rho$ changes from $-3h$ to $-h$
and $\xi$ from $-A$ to $A$ on $D_{2n-1}$ (see the figure). Notice that
this is different from the choice of coordinates we employed locally
in each cell in the previous section; in terms of the old coordinates,
$D_{2n-1}\subset C_{n-1}$ is bounded by the curves $\rho=h,$
$\rho=3h,$ $\xi=L-A,$ $\xi=L+A.$ With the new choice of coordinates,
we denote by $I_n$ the part of $C_{n-1} \cup C_{n}$ bounded by the
curves $\rho = -3h,$ $\rho =3h,$ $\xi=-A,$ and $\xi = A.$ We also
denote, in the extension of notation of the previous section,
\begin{equation}
\label{aver1}
 \langle T \rangle_{2n-1} = \frac{1}{2AF} \int\limits_{D_{2n-1}}
T(\rho, \xi) k(\rho)\, d\xi d\rho, 
\end{equation}
\begin{equation}
\label{aver2}
  \langle T \rangle_{2n} = \frac{1}{2AF} \int\limits_{D_{2n}}
T(\rho, \xi) k(\rho)\, d\xi d\rho 
\end{equation}
($k(\rho)$ is defined in a natural way as $k(-\rho)$ when $\rho<0$). 
 
\begin{thm}
\label{diff}
There exist constants $C_{1,2}$ depending only on the constants in the bounds
of Proposition \ref{procoord}, such that
\begin{equation}
\label{ediff}
\kappa \int\limits_{I_n} |\nabla T|^2 \, dxdy \geq \frac{C_1\kappa H}{h} 
\left( \langle T \rangle_{2n}- \langle T \rangle_{2n-1}\right)^2
\end{equation}
for $A=C_2H$.
\end{thm}

{\bf Proof.}
Notice that 
\begin{eqnarray*}
&&  \frac{1}{2AF} \left| \int\limits_{-A}^A d \xi 
\left( \int\limits_h^{3h} k(\rho) T(\rho, \xi) d \rho -
 \int\limits_{-h}^{-3h} k(\rho) T(\rho, \xi) d \rho \right) \right| 
 \leq  \frac{1}{2AF} \int\limits_{-A}^A d \xi 
\int\limits_h^{3h} k(\rho) d \rho  \int\limits_{-3h}^{3h} 
\left| \frac{\partial T}{\partial \rho}(\gamma, \xi) \right| d \gamma \\
&& \leq \left(\frac{3h}{A}\right)^{1/2}  
\left( \int\limits_{-A}^A  \int\limits_{-3h}^{3h} 
\left| \frac{\partial T}{\partial \rho}(\rho, \xi) \right|^2 
d \rho \right)^{1/2} \leq C \frac{h^{1/2}}{H^{1/2}} 
\left( \int\limits_{I_n} |\nabla T|^2 \, dxdy \right)^{1/2}. 
\end{eqnarray*}
Comparing the quantity we estimated with the definitions 
(\ref{aver1}) and (\ref{aver2}), we see that (\ref{ediff}) follows.
$\Box$

\section{Cellular flows: Reactive estimate}\label{sec:react}

It is worthy to note that the diffusive estimate we proved in the last
section is quadratic in the drop of the temperature over the interface. 
With the estimates we currently have we could not prove any lower 
bound for the burning rate, since the temperature could stay constant
inside the cells (so that advective estimate does not give us anything)
and drop in extremely small increments over the diffusive interfaces.
There is no lower bound for such scenario because of the quadratic 
dependence. But in this case, the region where $T \in [\theta_1, \theta_4]$ 
would be very large, and hence we can hope to get a lower bound on 
$\int f(T).$ The following theorem is a rigorous expression of the 
above idea. We are mostly interested in the case $H\ge\kappa/v_0$
but include the case $H\le {\frac{\kappa}{v_0}}$ for completeness.
\begin{thm}
\label{react}
Assume that in a given cell $C_n,$
\[ \frac{1}{F}
\int\limits_h^{3h} T(\rho, \xi) \, d \rho \in [\theta_1, \theta_4] \]
for every $\xi.$ 
Then 
\begin{equation}
\label{ereact}
\int\limits_{C_n} \left[ \kappa |\nabla T|^2 + \frac{v_0^2}{\kappa} 
f(T) \right] \, dxdy \geq C(f_0, \zeta) 
v_0 H {\rm min}(1, \frac{Hv_0}{\kappa}).
\end{equation}
\end{thm}
{\bf Proof.} The proof uses an argument we already used several 
times when proving advective estimate.
Consider a region in $C_n$ bounded by the curves $\rho=h$ and $\rho =H/2.$ 
By mean value theorem, for every $\xi,$ there exists $\rho_0 \in [h, 3h]$ 
where $T(\rho_0, \xi) \in [\theta_1, \theta_4].$ All $\xi$ fall into two
cases. In the first case, 
$T(\rho,\xi) \in [\theta_1-\zeta, \theta_4+\zeta]$ for 
every $\rho$ in the region $\rho\in[h,H/2]$. Then 
\[ \frac{v_0^2}{\kappa} \int\limits_h^{H/2} f(T(\rho, \xi)) d \rho
\geq \frac{f_0 H v_0^2}{\kappa}. \]
In the second case, there also exists $\rho_1 \in [h, H/2]$ 
such that $T(\rho_1, \xi)
\notin [\theta_1-\zeta, \theta_4+\zeta].$ In this case we can find
$\rho_2,$ $\rho_3$ such that $|T(\rho_2, \xi) - T(\rho_3, \xi)| \geq \zeta,$
and for every $\rho$ between $\rho_2$ and $\rho_3,$ we have 
$T(\rho, \xi) \in [\theta_1-\zeta, \theta_4+\zeta].$ In this case, by 
the usual argument, we have 
\[  f_0^{1/2}\zeta \leq \left( \int\limits_h^{H/2} f(T) \, d\rho
 \int\limits_h^{H/2} |\nabla T|^2  \, d\rho \right)^{1/2} \leq 
 \frac{1}{2v_0} \int\limits_h^{H/2} \left[ \kappa |\nabla T|^2 +
\frac{v_0^2}{\kappa} f(T) \right]\,d\rho. \]
It is easy to show following Lemma~\ref{Qest} and Proposition~\ref{procoord}
that the length of the interval of integration in $\xi,$
$2L,$ satisfies $2L \geq CH,$ where $C$ is some universal constant.  
Integrating over $\xi,$ we then get the result of the theorem.
$\Box$

\section{Cellular flows: Putting it together}\label{sec:putting}

Now we have all necessary estimates to establish the lower bound on
burning rate in the case of cellular flows. We begin with an
auxiliary computation. Its goal is to choose the right value of $h$
depending on other physical parameters fixed in the problem.  The
argument below may not be completely rigorous; we will fill in the
gaps in the actual proof.  Consider the sequence of averages
(\ref{aver1}), (\ref{aver2}) introduced in Section \ref{sec:diff}
devoted to the diffusive estimate:
\[\dots, \langle T \rangle_{2n},  \langle T \rangle_{2n+1}, 
 \langle T \rangle_{2(n+1)},  \langle T \rangle_{2n+3}, \dots \]
where $n$ varies from $-\infty$ to $\infty.$ The values of the 
averages tend to $1$ as $n \rightarrow -\infty$ and to $0$ as 
$n \rightarrow \infty.$ Assume that (for sufficiently large $U$)
the change is gradual, and there exists a number of consecutive 
numbers $n$ where 
all averages lie in the interval $[\theta_1, \theta_4].$ 
Moreover, assume that in all cells corresponding to these values 
of $n,$ we have 
\[ \frac{1}{F} \int_h^{3h} k(\rho) T(\rho, \xi)\, d\rho \in [\theta_1, 
\theta_4] \] 
for all $\xi,$ so that the reactive and linear advective 
estimates may be applied
in these cells. It does not concern us here that $h$ has to be chosen yet; 
as we mentioned above, this is an auxiliary computation and the rigorous
argument will appear in the proof. 
Let us denote
\[ \delta T_{n,a}= |\langle T \rangle_{2n+1}- \langle T \rangle_{2n}| \]
the change of the averages controlled by the advective estimate, and 
\[ \delta T_{n,d}= |\langle T \rangle_{2(n+1)}- \langle T \rangle_{2n+1}| \]
the drop of the temperature controlled by diffusive estimate.
We also set $\delta T_n = \delta T_{n,a}+\delta T_{n,d}$ which we call 
the total drop over $n$th cell. Then we have, according to the advective 
estimate (\ref{adv2}),
\begin{equation}
\label{advcontr}
\int\limits_{C_n} \left(T_t+\kappa |\nabla T|^2 + \frac{v_0^2}{\kappa}
f(T) \right) \, \frac{dxdy}{H} \geq C(f_0, \zeta) 
\left(1+\frac{\kappa}{v_0 h}\right)^{-1}
\frac{Uh}{H}  \delta T_{n,a}, 
\end{equation}
where the constant $C$ depends only on $f_0,$ $\zeta$ and fixed constants
associated with the geometry of streamlines of the flow. 
According to the diffusive estimate (\ref{ediff}), we also have 
\begin{equation}
\label{difcontr}
\int\limits_{C_n \cup C_{n+1}} 
\left(\kappa |\nabla T|^2 + \frac{v_0^2}{\kappa} f(T) \right) \, 
\frac{dxdy}{H} \geq C \frac{\kappa}{h} \delta T_{n,d}^2.
\end{equation}
Also for every $C_n$ from the region we consider, we have by 
the reactive estimate (\ref{ereact}) that
\begin{equation}
\label{reaccontr}
\int\limits_{C_n} \left(\kappa |\nabla T|^2 + \frac{v_0^2}{\kappa}
f(T) \right) \, \frac{dxdy}{H} \geq C(f_0, \zeta) v_0 {\rm min}(1, 
\frac{Hv_0}{\kappa}) = C(f_0, \zeta) v_0
\end{equation}
under the assumption $Hv_0/\kappa>1$ (see (\ref{eq:cellsize})).
Fix some  cell $C_n,$ and introduce a parameter
$0 \leq a_n \leq 1$ such that 
\[ 
a_n \delta T_n = \delta T_{n,d}, \,\,\, (1-a_n) \delta T_n =\delta T_{n,a}. \]
According to (\ref{advcontr}), (\ref{difcontr}), and (\ref{reaccontr}), we have
\begin{equation}
\label{total}
\int\limits_{C_n \cup C_{n+1}} 
\left(\kappa |\nabla T|^2 + \frac{v_0^2}{\kappa} f(T) \right) \, 
\frac{dxdy}{H} \geq  C\left[  (1+\frac{\kappa}{v_0 h})^{-1} \frac{Uh}{H}   
\delta T_n (1-a_n) +
\frac{\kappa a_n^2}{h} \delta T_{n}^2 + v_0 \right].
\end{equation}
Over our sequence of cells, the value of temperature averages falls by
a fixed amount, $\theta_4-\theta_1.$ We are going to assume that this
fall is gradual, and consider the lower bound on the contribution of
the cells $C_n$ and $C_{n+1}$ to the total burning rate normalized by
the temperature falloff $\delta T_n$ over the cell $C_n.$ Namely,
denote
\[ 
V_n =  \int\limits_{C_n \cup C_{n+1}} 
\left(T_t+\kappa |\nabla T|^2 + \frac{v_0^2}{\kappa} f(T) \right) \, 
\frac{dxdy}{H}, 
\]
then the burning rate $V$ satisfies
\[ 
V \geq \frac{1}{6} \sum\limits_n V_n. 
\]
Equation (\ref{total}) may be rewritten as
\begin{equation}\label{lowerwitha}
V_n \geq  C\left[  \left(1+\frac{\kappa}{v_0 h}\right)^{-1} 
\frac{Uh}{H}   (1-a_n) +
\frac{\kappa a_n^2}{h} 
\delta T_{n} + \frac{v_0}{\delta T_n}  \right] \delta T_n. 
\end{equation}
Moreover, 
\[ \sum \delta T_n \geq \theta_4-\theta_1, 
\]
where summation is taken over the set of cells where all estimates
apply.  In the lower bound (\ref{lowerwitha}), $a_n$ and $|\delta T_n|$
are out of our control, but we may choose $h$ optimally.  Hence if we
denote
\[  B(\kappa, v_0, U, H)=\max\limits_{0\leq h \leq H/6} 
\left\{\min\limits_{0\leq a_n \leq 1, 0\leq \delta T_n \leq 2} \left(  
(1+\frac{\kappa}{v_0 h})^{-1} \frac{Uh}{H}  (1-a_n) +
\frac{\kappa a_n^2}{h} \delta T_{n}+\frac{v_0}{\delta T_n}\right)\right\}, \]
then 
\[ V \geq CB(\kappa, v_0, U, H) (\theta_4-\theta_1). \]
Therefore our goal is to find $B(\kappa, v_0, U, H).$
The minimum in $\delta T_n$ is achieved
when 
\[ 
\delta T_n = {\rm min}(2,\frac{1}{a_n}\sqrt{\frac{v_0 h}{\kappa}}). 
\]
Let us consider two different regimes. \\
1. If $\frac{\kappa}{v_0 h} \leq 1$ then the minimum is achieved for $\delta T_n\geq1$ and
\[ 
B(\kappa, v_0, U, H) \geq C\max\limits_{0\leq h \leq H/3} 
\left\{ \min\limits_{0\leq a_n  \leq 1} \left(  
 \frac{Uh}{H}  (1-a_n) +
\frac{\kappa a_n^2}{h}  + v_0  \right) \right\}, 
\]
where $C$ is a universal constant. 
Choose $h$ out of the condition $Uh/H = \kappa/h,$ so that
\begin{equation}
\label{hchoice1}
h =\sqrt{\frac{\kappa H }{U}}. 
\end{equation}
Then evidently
\[ B(\kappa, v_0, U, H) \geq C \left(\sqrt{\frac{\kappa U}{H}}+v_0 \right)=
Cv_0 \left(\sqrt{\frac{\tau_c}{\tau_u}}+1\right) \]
with $\tau_c$ and $\tau_u$ defined in (\ref{eq:turnover}).
Given  (\ref{hchoice1}), condition $\frac{\kappa}{v_0 h} \leq 1$ translates into
\[ \frac{1}{v_0}\sqrt{\frac{\kappa U}{H}} =\frac{\tau_c}{\tau_u}\leq 1. \]
For this computation to apply, we also need to
ensure that $h\le H/6,$ that is,
\[ \sqrt{\frac{\kappa  }{U H}} \leq 1/6, \]
so that $h$ is in the acceptable range. This is a condition
 on the Peclet number. All our bounds will remain
valid if replace the choice of $h$ by $Ch$ with $C$ in some fixed
range; we just have to adjust the universal constants.
Therefore, for presentation purposes we will henceforth require 
the condition (\ref{eq:peclet}) 
\[ {\rm Pe} =\sqrt{\frac{UH}{\kappa}} \geq 1. \]
2. If  $\frac{\kappa}{v_0 h} \geq 1$ then, since $\delta T_n = {\rm
  min}(2,\frac{1}{a_n}\sqrt{\frac{v_0 h}{\kappa}}),$ we have that
\[ 
B(\kappa, v_0, U, H) \geq \max_{0\leq h \leq H/6}\min_{0\leq a_n\leq
  1} g(a_n),
\]
where 
\[ 
g(a_n)=\left\{\begin{matrix}{ 
\frac{Uh^2 v_0}{\kappa H}(1-a_n) + \frac{2\kappa a_n^2}{h}+\frac{v_0}{2}, & 
a_n \leq \frac{1}{2}\sqrt{\frac{v_0 h}{\kappa}}\cr 
\frac{Uh^2 v_0}{\kappa H}(1-a_n) + 2a_n \sqrt{\frac{v_0 \kappa}{h}},
& a_n \geq \frac{1}{2}\sqrt{\frac{v_0 h}{\kappa}}.} \end{matrix}\right.
\]
Note that
\[
\max_{0\leq h \leq H/6}\min_{0\leq a_n  \leq 1} 
\left( \frac{Uh^2 v_0}{\kappa H}(1-a_n) + 
\frac{2\kappa a_n^2}{h}+\frac{v_0}{2} \right) \geq C((Uv_0H^{-1}\kappa)^{1/3}+v_0)=
Cv_0\left( \left(\frac{\tau_c}{\tau_u}\right)^{1/3}+1 \right), \]
since we can choose $h$ from $Uh^2v_0/(\kappa H) = \kappa /h,$ which gives
\begin{equation}
\label{hchoice2}
h = \left( \frac{\kappa^2 H}{Uv_0} \right)^{1/3}. 
\end{equation}
Also
\[
\max_{0\leq h \leq H/3}\min_{0\leq a_n  \leq 1} \left( 
\frac{Uh^2 v_0}{\kappa H}(1-a_n) + 
2a_n \sqrt{\frac{v_0 \kappa}{h}}\right) \geq C(Uv_0^3 \kappa H^{-1})^{1/5}=
Cv_0 \left(\frac{\tau_c}{\tau_u}\right)^{1/5}, 
\]
since we can choose $h$ from $Uh^2v_0/(\kappa H) = \sqrt{v_0
  \kappa/h},$ which gives
\begin{equation}
\label{hchoice3}
h = \left( \frac{\kappa^3 H^2}{U^2 v_0} \right)^{1/5}. 
\end{equation}
Therefore,
\[
\max_{0\leq h \leq H/3}\min_{0\leq a_n  \leq 1} g(a_n) 
\geq Cv_0{\rm min} \left\{ \left(
    \left(\frac{\tau_c}{\tau_u}\right)^{1/3}+1 \right),
  \left(\frac{\tau_c}{\tau_u}\right)^{1/5}\right\}. 
\] 
Notice that
$\tau_c/\tau_u \geq 1$ in the regime we consider, since
$\tau_c/\tau_u$ is also equal to $(\kappa/v_0 h)^{3}$ with $h$ chosen
according to (\ref{hchoice2}) or $(\kappa/v_0 h)^{5/2}$ with $h$
chosen according to (\ref{hchoice3}).  Therefore, we get that the
second regime is characterized by
\[ \frac{\tau_c}{\tau_u} =\frac{1}{v_0}\sqrt{\frac{\kappa U}{H}} \geq 1, \]
and in this case
\[ B(\kappa, v_0, U, H) \geq Cv_0\left(\frac{\tau_c}{\tau_u}\right)^{1/5} .\]
The choice of $h$ is given by $(\ref{hchoice3}).$ In this regime, we
also need to assume that
\[ \left(\frac{\kappa^3  }{U^2 H^3 v_0}\right)^{1/5} \leq 1/6, \]
 so that $h$ is in the acceptable range. 
This, up to a constant, is a condition 
\[
\frac{l}{H}\le \hbox{Pe}^2
\]
which is satisfied provided (\ref{eq:cellsize}) and (\ref{eq:peclet})
hold.

Now we are finally ready to give a proof of Theorem~\ref{cellthm}. \\
{\bf Proof.}  Given the parameters $\kappa,$ $v_0,$ $H,$ and $U,$
define $h$ according to
\begin{eqnarray}
\label{hcho1}
h =\frac{1}{6}\sqrt{\frac{\kappa H }{U}}, & \tau_c/\tau_u \leq 1 \\
\label{hcho2}
 h = \frac{1}{6}\left( \frac{\kappa^3 H^2}{U^2 v_0} \right)^{1/5},
&  \tau_c/\tau_u \geq 1.
\end{eqnarray}
These choices ensure that $h \leq H/6$ provided that (\ref{eq:peclet}),
(\ref{eq:cellsize}) are satisfied. 
Consider the smallest integer number $m_1$ such that $\langle T
\rangle_i \leq \theta_4-(\theta_4-\theta_1)/10$ for all $i \geq m_1$.
We need to consider several options. The first is taken care
of by 
\begin{lemma}
\label{lfirst}
If $\langle T \rangle_{m_1} \leq \theta_4-(\theta_4-\theta_1)/5$  
then
the bounds of Theorem~\ref{cellthm} hold. Moreover, the same is true
if for some $m$, $m+1,$ we have $\langle T \rangle_{m},\langle T
\rangle_{m+1} \in [\theta_1, \theta_4]$ and
\[ 
|\langle T \rangle_{m+1}-\langle T \rangle_{m}| \geq (\theta_4-\theta_1)/10. 
\]   
\end{lemma}
{\bf Proof.}  Let us consider the case where $\langle T \rangle_{m_1}
\leq \theta_4-(\theta_4-\theta_1)/5$, the other case is similar. If
$m_1=2n$ is even, notice that by definition of $m_1$, $\langle T
\rangle_{2n-1} \geq \theta_4- (\theta_4-\theta_1)/10.$ Then by the
diffusive estimate (\ref{ediff}) we obtain
\[ 
V \geq \frac{1}{3}V_n \geq C \frac{\kappa}{h} (\theta_4-\theta_1)^2. 
\]
In the case where $\tau_c/\tau_u\leq 1$ putting in $h$ given by
(\ref{hcho1}) gives exactly the bound (\ref{perter1}) of
Theorem~\ref{cellthm}, modulo $C_2v_0.$ But by
Proposition~\ref{general}, we can always add $Cv_0$ to the lower bound
on $V$. In the case where $\tau_c/\tau_u\geq 1,$ putting in $h$ given by
(\ref{hcho2}) gives the lower bound $C(\kappa^2 U^2 H^{-2}
v_0)^{1/5}=Cv_0(\tau_c/\tau_u)^{2/5},$ which is even better than the
 bound (\ref{perter1}) (modulo $v_0$).  Hence,
Theorem \ref{cellthm} also holds.

If $m_1=2n+1$ is odd, then by definition of $m_1,$ $\langle T
\rangle_{2n} \geq \theta_4-(\theta_4-\theta_1)/10$. 
Thus temperature falls inside the
cell $C_n$, and we may use advective estimate of Theorem~\ref{1adv}.
We can find
$\xi_0,$ $\xi_1$ in the cell $C_n$ such that
\[  
\frac{1}{F} \int\limits_h^{3h} k(\rho) T(\rho, \xi_0)\,d \rho =\theta_4-
\frac{\theta_4-\theta_1}{10},
 \,\,\,
\frac{1}{F} \int\limits_h^{3h} k(\rho) T(\rho, \xi_1)\,d \rho =
\theta_4-\frac{\theta_4-\theta_1}{5}, \]
and 
\[ 
\frac{1}{F} \int\limits_h^{3h} k(\rho) T(\rho, \xi)\,d \rho 
\in [\theta_4-\frac{\theta_4-\theta_1}{5}, \theta_4-
\frac{\theta_4-\theta_1}{10}]
 \]
for every $\xi$ between $\xi_0$ and $\xi_1.$ By the advective estimate
(\ref{advest1}), we get
\[ 
V\geq \frac{1}{3}
V_n \geq C(\zeta, f_0) \left(1+ \frac{\kappa}{v_0 h} \right)^{-1} 
(\theta_4 -\theta_1)^3 \frac{Uh}{H}. 
\]
Direct substitution of the expression (\ref{hcho1}) or
(\ref{hcho2}) for $h$ depending on the value of
$\tau_c/\tau_u$ and comparison of the above bound with (\ref{perter1})
gives the conclusion of Theorem \ref{cellthm}.
The proof of the second statement of this lemma is parallel to the
above argument.  $\Box$

Lemma~\ref{lfirst} proves Theorem~\ref{cellthm} unless there exists a
sequence $m_1, \dots, m_2$ (of length at least $8,$ in fact) such
that for every $m_1 \leq m \leq m_2,$
\[ 
\langle T \rangle_{m} \in [\theta_1+\frac{\theta_4-\theta_1}{10}, 
\theta_4-\frac{\theta_4-\theta_1}{10}], 
\]
and for $m=m_1, \dots, m_2-1$
\[ 
|\langle T \rangle_{m+1} - \langle T \rangle_{m}| \leq \frac{\theta_4-\theta_1}{10}.
\] 
Therefore there exists a sequence of at least three (or more)
consecutive cells $C_{n_1}, \dots, C_{n_2}$ such that
\begin{equation}
\label{vary}
\langle T \rangle_{2n_1} \geq \theta_4 - \frac{3}{10}(\theta_4-\theta_1),
\,\,\, \langle T \rangle_{2n_2+1} \leq
 \theta_1+\frac{3}{10}(\theta_4-\theta_1), 
\end{equation}
and 
\begin{equation}
\label{catch}
\langle T \rangle_{m} \in [\theta_4-\frac{\theta_4-\theta_1}{10}, 
\theta_1+\frac{\theta_4-\theta_1}{10}]
\end{equation}
for every $m$ such that $2n_1 \leq m \leq 2n_2+1$. Assume that there
exists a cell $C_n$ with $n_1\leq n \leq n_2$ where we can find
$\xi_0$ so that
\[ 
\int\limits_h^{3h} k(\rho) T(\rho,  \xi_0)\, 
d\rho \notin [\theta_1, \theta_4]. 
\]
By (\ref{catch}) and mean value theorem we can also find $\xi_1$ in
this cell such that
\[ 
\int\limits_h^{3h} k(\rho) T(\rho,  \xi_1)\, d\rho \in 
[\theta_1+\frac{\theta_4-\theta_1}{10}, \theta_4
-\frac{\theta_4-\theta_1}{10}]. 
\] 
Therefore, the advective estimate (\ref{advest1}) can be applied in
this cell, giving
\[ 
V_n \geq   C(\zeta, f_0) \left(1+ \frac{\kappa}{v_0 h} 
\right)^{-1} (\theta_4 -\theta_1)^3 \frac{Uh}{H}.
\]
Hence, similarly to the proof of Lemma~\ref{lfirst},
Theorem~\ref{cellthm} holds in this case.  The only case left to
consider is the case where for every $n_1 \leq n \leq n_2,$ for every
$\xi$ in a cell $C_n,$ we have
\[ 
\int\limits_h^{3h} k(\rho) T(\rho,  \xi)\, d\rho \in [\theta_1, \theta_4]. 
\]
In this case the second advective estimate (\ref{adv2}), as well as
reactive estimate (\ref{ereact}) apply in every cell $C_n$ such that
$n_1 \leq n \leq n_2.$ Recall the notation
\[ 
(1-a_n) \delta T_n =\delta T_{n,a} =  |\langle T \rangle_{2n+1} -\langle T 
\rangle_{2n}|, \,\,\, 
a_n \delta T_n =\delta T_{n,d} =  |\langle T \rangle_{2(n+1)} -\langle T 
\rangle_{2n+1}|. 
\]
Following the computation we performed at the beginning of this
section we get
\begin{equation}
\label{kill}
V \geq \frac{1}{6}
\sum\limits_{n=n_1}^{n_2} V_n \geq C \sum\limits_{n=n_1}^{n_2}
\left[ \left( \left( (1+\frac{\kappa}{v_0 h} \right)^{-1} \frac{Uh}{H}(1-a_n)
+ \frac{\kappa}{h} a_n^2 \delta T_n 
+v_0 \delta T_n^{-1} \right) \delta T_n \right]. 
\end{equation}
Consider the case where $\tau_c/\tau_u\geq 1$ (the other case
is similar but simpler). Putting the expression for $h$ from (\ref{hcho2})
into (\ref{kill}), we get 
\[ 
V \geq  C \sum\limits_{n=n_1}^{n_2}
\left( \left[  v_0\left(\frac{\tau_c}{\tau_u}\right)^{1/5}(1-a_n)+ 
v_0\left( \frac{\tau_c}{\tau_u}\right)^{2/5} a_n^2 \delta T_n 
+v_0 \delta T_n^{-1} \right] \delta T_n \right). 
\]
Since by (\ref{vary}), $\sum_{n=n_1}^{n_2} \delta T_n
\geq 2(\theta_4-\theta_1)/5,$ we have
\begin{equation}
\label{kill1}
 V \geq  C v_0\min\limits_{0\leq a_n\leq 1, 0\leq \delta T_n  \leq 2}
 \left[  \left(\frac{\tau_c}{ \tau_u}\right)^{1/5}(1-a_n)+ 
\left( \frac{\tau_c}{\tau_u}\right)^{2/5} a_n^2 \delta T_n 
+\delta T_n^{-1} \right]. 
\end{equation}
It remains to show that the expression 
in square brackets in (\ref{kill1}) is always greater or equal to 
$Cv_0(\tau_c/\tau_u)^{1/5}$ (we can always add $Cv_0$ later to the lower 
bound for $V$ by Proposition~\ref{general}). 
If $a_n<1/2,$ the first term in the sum gives
exactly the estimate we need (no matter what is the value of 
$\delta T_n \leq 1$). Hence it remains to consider the 
case where $a_n>1/2.$ 
In this case the sum of the second and third term in the square brackets
is greater than or equal to
\[ 
Cv_0\left( \left( \frac{\tau_c}{\tau_u}\right)^{2/5} \delta T_n 
+\delta T_n^{-1} \right) \geq Cv_0\left(\frac{\tau_c}{\tau_u}\right)^{1/5} 
\]
and this finishes the proof of Theorem \ref{cellthm}. $\Box$

\end{document}